\documentclass[final,3p,times,fleqn]{elsarticle}
\usepackage{amsmath,amssymb,amscd,mathtools}
\usepackage{enumitem}
\usepackage{graphicx}
\usepackage{xcolor}
\usepackage{algorithm}
\usepackage{algpseudocode}
\usepackage{diagbox}
\usepackage{url}
\usepackage{tikz}
\usepackage{xparse}
\usepackage{booktabs}

\usepackage[colorlinks=true, allcolors=blue]{hyperref}
\usepackage{hypernat}

\usepackage{subcaption}
\usepackage{pgfplots}
\usepackage{pgfplotstable}
\definecolor{butter1}{rgb}{0.988,0.914,0.310}
\definecolor{chocolate1}{rgb}{0.914,0.725,0.431}
\definecolor{chameleon1}{rgb}{0.541,0.886,0.204}
\definecolor{skyblue1}{rgb}{0.247,0.524,0.912}
\definecolor{applegreen}{rgb}{0.55, 0.71, 0.0}
\definecolor{blue-green}{rgb}{0.0, 0.87, 0.87}
\definecolor{plum1}{rgb}{0.678,0.498,0.659}
\definecolor{scarletred1}{rgb}{0.937,0.161,0.161}

\pgfplotsset{
        colormap={jet-like-GLVis}{
rgb255=(0,0,127.5)
rgb255=(0,0,143.4375)
rgb255=(0,0,159.375)
rgb255=(0,0,175.3125)
rgb255=(0,0,191.25)
rgb255=(0,0,207.1875)
rgb255=(0,0,223.125)
rgb255=(0,0,239.0625)
rgb255=(0,0,255)
rgb255=(0,15.9375,255)
rgb255=(0,31.875,255)
rgb255=(0,47.8125,255)
rgb255=(0,63.75,255)
rgb255=(0,79.6875,255)
rgb255=(0,95.625,255)
rgb255=(0,111.5625,255)
rgb255=(0,127.5,255)
rgb255=(0,143.4375,255)
rgb255=(0,159.375,255)
rgb255=(0,175.3125,255)
rgb255=(0,191.25,255)
rgb255=(0,207.1875,255)
rgb255=(0,223.125,255)
rgb255=(0,239.0625,255)
rgb255=(0,255,255)
rgb255=(0,248.625,223.125)
rgb255=(0,242.25,191.25)
rgb255=(0,235.875,159.375)
rgb255=(0,229.5,127.5)
rgb255=(0,223.125,95.625)
rgb255=(0,216.75,63.75)
rgb255=(0,210.375,31.875)
rgb255=(0,204,0)
rgb255=(73.236,210.375,0)
rgb255=(111.0015,216.75,0)
rgb255=(141.576,223.125,0)
rgb255=(168.249,229.5,0)
rgb255=(192.3465,235.875,0)
rgb255=(214.5825,242.25,0)
rgb255=(235.365,248.625,0)
rgb255=(255,255,0)
rgb255=(255,239.0625,0)
rgb255=(255,223.125,0)
rgb255=(255,207.1875,0)
rgb255=(255,191.25,0)
rgb255=(255,175.3125,0)
rgb255=(255,159.375,0)
rgb255=(255,143.4375,0)
rgb255=(255,127.5,0)
rgb255=(255,111.5625,0)
rgb255=(255,95.625,0)
rgb255=(255,79.6875,0)
rgb255=(255,63.75,0)
rgb255=(255,47.8125,0)
rgb255=(255,31.875,0)
rgb255=(255,15.9375,0)
rgb255=(255,0,0)
rgb255=(239.0625,0,0)
rgb255=(223.125,0,0)
rgb255=(207.1875,0,0)
rgb255=(191.25,0,0)
rgb255=(175.3125,0,0)
rgb255=(159.375,0,0)
rgb255=(143.4375,0,0)
rgb255=(127.5,0,0)  
        },
    } 
    
 \pgfplotsset{
        colormap={calewhite-GLVis}{   
        rgb255=(0,0,255)
rgb255=(0,255,255)
rgb255=(0,255,0)
rgb255=(255,255,0)
rgb255=(255,127.5,0)
rgb255=(255,0,0)
rgb255=(255,0,255)
rgb255=(127.5,0,255)
        },
    } 

\def\generateTikzFigures{0}
\newcounter{FigureCounter}
\graphicspath{ {./tikz/} }

\usetikzlibrary{matrix,backgrounds,decorations.pathreplacing}
\pgfdeclarelayer{myback}
\pgfsetlayers{myback,background,main}

\tikzset{mycolor/.style = {line width=1bp,color=#1}}%
\tikzset{myfillcolor/.style = {draw,color=#1,fill=#1}}%

\NewDocumentCommand{\highlight}{O{blue!40} m m}{%
\draw[mycolor=#1] (#2.north west)rectangle (#3.south east);
}

\NewDocumentCommand{\fhighlight}{O{blue!40} m m}{%
\draw[myfillcolor=#1] (#2.north west)rectangle (#3.south east);
}

\biboptions{sort&compress}
\def\grid{{\mathcal{T}}}
\def\faces{{\mathcal{E}}}

\def\blkInterpolate{{\blkMat{P}_{\ell+1}^\ell}}
\def\blkRestrict{{\blkMat{R}_\ell^{\ell+1}}}

\def\blkInject{{\blkMat{Q}_\ell^{\ell+1}}}

\newcommand{\interpolate}[1]{ \big(\Mat{P}_{#1}\big)_{\ell+1}^\ell }
\newcommand{\restrict}[1]{ \big(\Mat{R}_{#1}\big)_\ell^{\ell+1}}

\newcommand{\inject}[1]{ \big(\Mat{Q}_{#1}\big)_\ell^{\ell+1}}

\newcommand{\interpolateOne}[1]{ \Mat{P}_{#1} }
\newcommand{\restrictOne}[1]{ \Mat{R}_{#1} }

\newcommand{\diag}[1]{{\bf diag }\left({#1}\right)}

\def\TimeDomain{{\mathbb{T}}}

\newcommand{\tensorOne}[1]{\boldsymbol{#1}}

\renewcommand{\Vec}[1]{%
  \ifcat\noexpand#1\relax 
    \boldsymbol{#1}
  \else
    \mathbf{#1}
  \fi
}
\newcommand{\Mat}[1]{#1}

\newcommand{\blkVec}[1]{\Vec{#1}}
\newcommand{\blkMat}[1]{\Vec{#1}}

\newdefinition{rmk}{Remark}
\newdefinition{definition}{Definition}
\usepackage{amsthm}
        \theoremstyle{plain}
        \newtheorem{thm}{Theorem}
        \newtheorem{proposition}[thm]{Proposition}
       	\newtheorem{lemma}{Lemma}

\newcommand{\chak}[1]{{\color{black}#1}}

\begin{document}
\title{An aggregation-based nonlinear multigrid solver for two-phase flow and transport in porous media}
\begin{frontmatter}
\author[casc]{Chak Shing Lee\corref{cor1}}
\ead{cslee@llnl.gov}
\author[total]{Fran\c cois P. Hamon}
\ead{francois.hamon@totalenergies.com}
\author[aeed]{Nicola Castelletto}
\ead{castelletto1@llnl.gov}
\author[casc,psu]{Panayot S. Vassilevski}
\ead{vassilevski1@llnl.gov, panayot@pdx.edu}
\author[aeed]{Joshua A. White}
\ead{white230@llnl.gov}

\cortext[cor1]{Corresponding author.}
\address[casc]{Center for Applied Scientific Computing, Lawrence Livermore National Laboratory, Livermore, CA 94550, USA}
\address[total]{TotalEnergies E\&P Research and Technology, Houston, TX 77002, USA}
\address[aeed]{Atmospheric, Earth, and Energy Division, Lawrence Livermore National Laboratory, Livermore, CA 94550, USA}
\address[psu]{Fariborz Maseeh Department of Mathematics and Statistics, Portland State University, Portland, OR 97201, USA} 

\begin{abstract}
A nonlinear multigrid solver for two-phase flow and transport in a mixed fractional-flow velocity-pressure-saturation formulation is proposed.
The solver, which is under the framework of the full approximation scheme (FAS), extends our previous work on nonlinear multigrid for heterogeneous diffusion problems.
The coarse spaces in the multigrid hierarchy are constructed by first aggregating degrees of freedom, and then solving some local flow problems.
The mixed formulation and the choice of coarse spaces allow us to assemble the coarse problems without visiting finer levels during the solving phase, which is crucial for the scalability of multigrid methods.
Specifically, a natural generalization of the upwind flux can be evaluated directly on coarse levels using the precomputed coarse flux basis vectors.
The resulting solver is applicable to problems discretized on general unstructured grids.
The performance of the proposed nonlinear multigrid solver in comparison with the standard single level Newton's method is demonstrated through challenging numerical examples.
It is observed that the proposed solver is robust for highly nonlinear problems and clearly outperforms Newton's method in the case of high Courant-Friedrichs-Lewy (CFL) numbers.
\end{abstract}


\end{frontmatter}


\allowdisplaybreaks{}

\section{Introduction}

Numerical simulation of subsurface flow and transport is important for applications like petroleum recovery and CO$_2$ sequestration.
The problem is challenging due to the high degree of nonlinearity in the governing equations, the complex geometry of geological features, and the heterogeneity in rock properties.
This high heterogeneity yields a range of CFL numbers that can span orders of magnitude across the computational domain.
%
%
%
We consider here the fully implicit method (FIM), in which all the unknowns are treated implicitly in the time-stepping scheme.
In this case, a fully coupled system of discrete nonlinear equations needs to be solved at every time step.
%
%
In realistic field-scale simulations, solving these large, ill-conditioned systems is challenging and computationally expensive, especially when nonlinear convergence is slow.
Therefore, scalable and robust solvers that can be applied to a broad class of discretizations---on both structured and unstructured grids---are very desirable.

Over the past few decades, significant effort has been devoted to developing nonlinear solvers for the systems arising from FIM.
A popular approach is to first linearize the system using Newton's method or related variants, and then to apply a scalable linear solver to the Jacobian system \cite{cao2005}.
For example, if a linear multigrid solver is used to solve the Jacobian system, the overall method is referred to as Newton-multigrid \cite{henson2003multigrid}.
However, the performance of Newton-based solution algorithms can by crippled by slow nonlinear convergence.
In particular, the convergence of Newton's method can deteriorate significantly for poor initial guesses, which is problematic when large time step sizes are selected.
As a result, globalization methods \cite{deuflhard2011newton} have been developed to improve the robustness of Newton's method and avoid convergence failures.
Damping strategies for the Newton updates have been used extensively to enlarge the convergence radius.
They include local saturation chopping strategies based on heuristics \cite{younis2011modern}, as well as physics-based damping methods that can, in some cases, achieve unconditional convergence, see, e.g., \cite{jenny2009, wang2013, li2015} and \cite{moyner2017}.
%

In recent years, alternative nonlinear strategies have been applied to multiphase flow problems to overcome the limitations of Newton-based methods.
They include homotopy continuation methods \cite{younis2009, jiang2018}, in which robustness is achieved by solving an easier problem that is gradually relaxed towards the original problem.
%
%
These approaches can handle very large time step sizes and prevent convergence failures.
In ordering-based methods \cite{kwok2007, natvig2008, hamon2016ordering, kelmetstal2020reordering}, the degrees of freedom are reordered based on the phase potential direction to obtain a (block) triangular transport nonlinear system. 
These blocks are then solved sequentially, one at a time, which greatly accelerates nonlinear convergence and reduces computational cost.
Nonlinear preconditioning based on additive and multiplicative Schwarz preconditioned inexact Newton (A/MSPIN) \cite{cai2002, liu2015, dolean2016} is another relevant class of methods that has successfully been applied to multiphase flow and transport \cite{skogestad2013, skogestad2016, kelmetstal2020schwarz}.
Finally, nonlinear multigrid, which directly applies the multigrid concept at the nonlinear level (unlike Newton-multigrid), is an attractive alternative as a scalable solver.
Some nonlinear multigrid solvers based on full approximation scheme \cite{brandt77} for porous media flow were proposed in \cite{christensen16, christensen18, toft18, fas-spectral-diffusion}.
In particular, it was observed in \cite{fas-spectral-diffusion} that nonlinear multigrid can outperform Newton-multigrid as the underlying problem becomes stiffer.

In this paper, our goal is to develop a nonlinear multigrid solver for two-phase flow and transport problems in the subsurface,
building upon our previous work on nonlinear multigrid for heterogeneous diffusion problems \cite{fas-spectral-diffusion}.
To this end, a mixed fractional-flow velocity-pressure-saturation formulation of the two-phase flow and transport problem is considered,
where the primary unknowns in the discrete system are the total flux, pressure and saturation (of the wetting phase).
For the coarsening of the total flux and pressure, we adopt the lowest order version of the coarsening method in \cite{fas-spectral-diffusion}.
As for the saturation, coarse spaces are composed of piecewise-constant functions on algebraically constructed coarse grids.
The formulation and the coarse spaces allow us to assemble the coarse problems directly on the coarse levels using precomputed quantities.
Specifically, the formulation naturally leads to a generalization of the usual upwind direction selection operator on coarse levels.
This in turn enables us to derive the Jacobian system in a compact form.
We show, using challenging benchmark problems, that the proposed multigrid solver exhibits a more robust nonlinear convergence behavior than Newton's method, especially for large time steps.
The improved robustness, combined with the scalability of the multigrid methodology, results in significant reductions in the computational cost of the simulations.

The remainder of the paper is organized as follows.
In Section~\ref{sec:model_problem}, the system of nonlinear partial differential equations (PDEs) of interest and its finite volume discretization are described.
Then, the components of the proposed nonlinear multigrid solver are discussed in details in Section~\ref{sec:multigrid}.
Numerical examples comparing the performance of the proposed nonlinear multigrid to Newton's method are presented in Section~\ref{sec:numerics}.
Lastly, some conclusions are drawn in Section~\ref{sec:conclusion}.

\section{Model problem} \label{sec:model_problem}

We consider a two-phase flow and transport problem involving two immiscible and incompressible phases---a wetting phase, $w$, and a non-wetting phase, $nw$---flowing in an incompressible porous medium.
We focus on a mixed fractional-flow velocity-pressure-saturation formulation.
We neglect gravitational and capillary forces, a frequent assumption in many practical engineering applications.
Therefore, the pressure is the same for both phases, i.e. $p_w = p_{nw} = p$.
In this work, using the saturation constraint $\sum_{\alpha=\{w,nw\}} s_{\alpha} = 1$, we use the wetting-phase saturation as primary unknown and denote it from now on as $s = s_w$.

For a simply-connected polygonal domain  $\Omega \in \mathbb{R}^3$ and time interval $\TimeDomain := (T_0, T_f)$, with $T_0$ and $T_f$ the initial and final time, respectively, the strong form of the initial/boundary value problem (IBVP) consists of finding the total Darcy velocity $\tensorOne{v}: \Omega \times \TimeDomain \rightarrow \mathbb{R}^3$, the pressure $p: \Omega \times \TimeDomain \rightarrow \mathbb{R}$, and the wetting-phase saturation $s: \Omega \times \TimeDomain \rightarrow \mathbb{R}$ such that \cite{aziz79}:
%
\begin{subequations}
\begin{align}
	&\frac{1}{\lambda (s)} \mathbb{K}^{-1} \cdot \tensorOne{v} + \nabla p
	= 0
	&& \text{in} \; \Omega \times \TimeDomain
	&& \mbox{(total Darcy velocity)} ,
	\label{eq:total_darcy} \\
	&\nabla \cdot \tensorOne{v}
	=
	q^I(p, s) - q^P(p, s)
	&& \text{in} \; \Omega \times \TimeDomain
	&& \mbox{(total volume conservation)},
	\label{eq:total_mass_conservation} \\
	&\phi \frac{\partial s}{\partial t} + 
  \nabla\cdot [ f_w (s) \tensorOne{v} ]
	=
	q_w^I(p,s) - q_w^P(p,s),
	&& \text{in} \; \Omega \times \TimeDomain
	&& \mbox{(wetting-phase volume conservation)} ,
	\label{eq:wetting_mass_conservation}
\end{align}
\label{eq:model}\null
\end{subequations}
%
where
\begin{itemize}
  \item $\lambda(s) = \sum_{\alpha=\{w,nw\}} \lambda_{\alpha}(s)$ is the total mobility, with the corresponding phase-based quantities defined as the ratio of relative permeability, $k_{r,\alpha}$, to viscosity, $\mu_\alpha$, i.e. $\lambda_\alpha(s) := k_{r,\alpha}(s) / \mu_\alpha$. Various constitutive relationships for $k_{r,\alpha}$ will be considered in our numerical examples. Note that using standard assumptions on the phase mobilities, the total mobility is bounded away from zero;
  \item $\mathbb{K}$ and $\phi$ are the medium absolute permeability tensor and porosity, respectively;
  \item $q^I(p,s) = \sum_{\alpha=\{w,nw\}} q^I_{\alpha}(p,s)$ is the total volumetric source per unit volume, with $q^I_{\alpha}(p,s)$ the corresponding phased-based quantity. The term $q^P(p,s)$ is defined similarly and represents a total volumetric sink term. In this work, such terms are introduced to model wells based on inflow-performance relationships that depend on $p$ and $s$ \cite{lie19};
  \item $f_w(s) := \lambda_w(s) / \lambda(s)$ is the fractional flow function.
\end{itemize}
Without loss of generality, in our simulations the domain boundary, $\partial \Omega$, is always subject to no-flow boundary conditions.
This represents a natural assumption when simulating closed-flow systems, e.g. reservoirs containing petroleum fluids.
To ensure uniqueness of the pressure solution, we prescribe a datum value for pressure internally in the domain through sink and/or source terms.
The formulation is completed by appropriate initial conditions for $\tensorOne{v}$, $p$, and $s$.

\begin{rmk}
In the reservoir simulation community, Eqs. \eqref{eq:total_darcy}-\eqref{eq:total_mass_conservation} are the mixed form of what is typically referred to as the \textit{pressure equation}, whereas Eq. \eqref{eq:wetting_mass_conservation} is often called the \textit{saturation equation}.
\end{rmk}

\subsection{Finite Volume discretization}\label{sec:FV_tpfa}
The system of PDEs \eqref{eq:model} is discretized by a cell-centered two-point flux approximation (TPFA) finite-volume (FV) method \cite{EymGalHer00} on a conforming triangulation of the domain, combined with the backward Euler (fully implicit) time-stepping scheme.
First, we introduce some notation.
Let $\mathcal{T}$ be the set of cells in the computational mesh such that $\overline{\Omega} = \sum_{\tau \in \mathcal{T}} \overline{\tau}$.
For a cell $\tau_K \in \mathcal{T}$, with $K$ a global index, let $|\tau_K|$ denote the volume, $\partial \tau_{K} = \overline{\tau}_{K} \setminus \tau_{K}$ the boundary, $\tensorOne{x}_{K}$ the barycenter, and $\tensorOne{n}_{K}$ the outer unit normal vector associated with $\tau_K$.
Let $\mathcal{E}$ be the set of internal faces in the computational mesh included in $\Omega$.
An internal face $\varepsilon$ shared by cells $\tau_K$ and $\tau_L$ is denoted as $\varepsilon_{K,L} = \partial \tau_K \cap \partial \tau_L$, with the indices $K$ and $L$ such that $K < L$.
%
%
The area of a face is $|\varepsilon|$.
A unit vector $\tensorOne{n}_{\varepsilon}$ is introduced to define a unique orientation for every face, and we set $\tensorOne{n}_{\varepsilon} = \tensorOne{n}_{K}$.
To indicate the mean value of a quantity $(\cdot)$ over a face $\varepsilon$ or a cell $\tau_K$, we use the notation $(\cdot)_{\left| \right. \varepsilon}$ and $(\cdot)_{\left| \right. K}$, respectively.
Let $T_0 = t_0 < t_1 < \cdots < t_n = T_f$ be a partition of the time domain $\TimeDomain$.
The discrete (finite-difference) approximation to a time-dependent quantity $\chi(t_m)$ at time $t_m$ is denoted by $\chi^m$.
Also, we define the time step size $\Delta t_m := t_m - t_{m-1}$.

We consider a piecewise-constant approximation for both pressure and saturation.
For each cell $\tau_{K} \in \mathcal{T}$, we introduce one pressure, $p_K$, and one saturation, $s_K$, degree of freedom, respectively.
We denote by $\sigma_{\varepsilon}$ the numerical flux approximating the total Darcy flux through an internal face $\varepsilon=\varepsilon_{K,L}$, i.e. $\sigma_{\varepsilon} \approx \int_{\varepsilon_{K,L}} \tensorOne{v} \cdot \tensorOne{n}_{\varepsilon} \mathrm{d}\Gamma$, such that:
%
\begin{equation}
  \left( \frac{1}{\lambda(s_K) \overline{\Upsilon}_{K,\varepsilon}} + \frac{1}{\lambda(s_L) \overline{\Upsilon}_{L,\varepsilon}} \right) \sigma_{\varepsilon} \chak{-  (p_{K} - p_{L})} = 0,  
  \label{eq:TPFA_flux}
\end{equation}
%
where $\overline{\Upsilon}_{K,\varepsilon} (s_K)$ and $\overline{\Upsilon}_{L,\varepsilon} (s_L)$ are the constant (geometric) one-sided transmissibility coefficients, defined as \cite{lie19}
%
\begin{align}
  \overline{\Upsilon}_{i,\varepsilon} &=
  |\varepsilon_{K,L}| \frac{\tensorOne{n}_{i} \cdot \mathbb{K}_{\left| \right. i} \cdot (\tensorOne{x}_{\varepsilon} - \tensorOne{x}_{i})}{||\tensorOne{x}_{\varepsilon} - \tensorOne{x}_{i}||_2^2},
  &
  i &= \{K, L \},
  \label{eq:half-transmissibility}
\end{align}
%
with $\tensorOne{x}_{\varepsilon}$ a collocation point introduced for every $\varepsilon \in \mathcal{E}$ to enforce point-wise pressure continuity across interfaces.

\noindent
The approximation of the wetting-phase Darcy flux through $\varepsilon = \varepsilon_{K,L} \in \mathcal{E}$ in the discrete form of Eq.~\eqref{eq:wetting_mass_conservation} relies on using single-point upstream weighting (SPU) according to the sign of $\sigma_{\varepsilon}$, namely
%
\begin{align}
  f_w^{\text{upw}}(s_K,s_L) \sigma_{\varepsilon} 
  &\approx \int_{\varepsilon_{K,L}} f_w(s)\tensorOne{v}\cdot\tensorOne{n}_{\varepsilon}\mathrm{d}\Gamma,
  &
  f_w^{\text{upw}}(s_K,s_L)
  &=
  \begin{cases}
    f_w(s_K), & \text{if } \sigma_{\varepsilon} > 0, \\
    f_w(s_L), & \text{otherwise}.
  \end{cases}
  \label{eq:upwinding}
\end{align}
%
Source and sink terms in Eqs.~\eqref{eq:total_mass_conservation}-\eqref{eq:wetting_mass_conservation} are used to simulate the effect of injection and production wells.
We employ a conventional Peaceman well model \cite{Pea78}, which relates well control parameters, such as bottomhole pressure (BHP), to flow rates through the wellbore \cite{Pea78}.
We assume each well segment to be vertical, with a single perforation connected to the centroid of a cell.
Also, without lost of generality, we restrict ourselves to rate-controlled injection wells and BHP-controlled production wells.
For a cell $\tau_K$ connected to a well, source/sink terms are expressed as
%
\begin{align}
  q_{\alpha}^I(p_K, s_K) &=
    \bar{q}_{\alpha}^I \delta ( \tensorOne{x} - \tensorOne{x}_K ),
  \\  
  q_{\alpha}^P(p_K, s_K) &=
  - \lambda_{\alpha}(s_K) WI (\bar{p}_{bh} - p_K) \delta ( \tensorOne{x} - \tensorOne{x}_K ),
  \label{eq:peaceman}
\end{align}
%
where $\bar{q}_{\alpha}$ is the known $\alpha$-phase rate control, $\delta ( \tensorOne{x} - \tensorOne{x}_K )$ is the Dirac function, $WI$ is the well Peaceman index, and $\bar{p}_{bh}$ is the prescribed bottomhole pressure.
A comprehensive presentation on well models and well index calculation can be found in \cite{CheHuaMa06}.
For each cell $\tau_K \in \mathcal{T}_P$, with $\mathcal{T}_P$ the set of cells connected to a production well, we define the integral total volumetric production flux $\sigma_{K}^P = \int_{\tau_K} ( q^P_w + q^P_{nw} )\mathrm{d}\Omega$ such that
%
\begin{equation}
  \frac{1}{\lambda(s_K) WI} \sigma^P_{K} \chak{-  (p_{K} - \bar{p}_{bh} )} = 0.
  \label{eq:production_well_flux}
\end{equation}
%
Introducing coefficient vectors
$\Vec{\sigma}^m = \begin{bmatrix} \Vec{\sigma}_{\varepsilon}^m \\ \Vec{\sigma}_{\tau}^{m} \end{bmatrix}$,
$\Vec{\sigma}_{\varepsilon}^m = (\sigma_{\varepsilon}^m)_{\varepsilon \in \mathcal{E}}$,
$\Vec{\sigma}_{\tau}^m = (\sigma_{K}^{P,m})_{\tau_K \in \mathcal{T}_P}$,
$\Vec{p}^m = (p_{K}^m)_{\tau_{K} \in \mathcal{T}}$ and
$\Vec{s}^m = (s_{K}^m)_{\tau_{K} \in \mathcal{T}}$
that contain the unknown degrees of freedom at time $t = t_m$ (i.e. face fluxes, production well fluxes, cell pressures, and cell saturations) the matrix form associated with the IBVP \eqref{eq:model} can be stated as follows: 
given the discrete solution $\Vec{x}^{m-1} = \{ \Vec{\sigma}^{m-1}, \Vec{p}^{m-1}, \Vec{s}^{m-1} \}$ at time $t = t_{m-1}$, find $\Vec{x}^{m} = \{ \Vec{\sigma}^{m}, \Vec{p}^{m}, \Vec{s}^{m} \}$ such that
%
\begin{equation}
\Vec{r}^m(\Vec{x}^{m}) 
:=
\begin{bmatrix}
\Vec{r}_\sigma^m(\Vec{\sigma}^{m}, \Vec{p}^{m}, \Vec{s}^{m}) \\
\Vec{r}_p^m(\Vec{\sigma}^{m}) \\
\Vec{r}_s^m(\Vec{\sigma}^{m}, \Vec{s}^{m})  
\end{bmatrix}
: =
\begin{bmatrix}
M(\Vec{s}^{m})\Vec{\sigma}^{m} \chak{- D^T} \Vec{p}^{m} - \Vec{g}^{m}\\
D\Vec{\sigma}^{m} - \Vec{f}^{m} \\
T^m(\Vec{\sigma}^{m}, \Vec{s}^{m}) - (\Delta t_m)^{-1}W\Vec{s}^{m-1} - \Vec{h}^{m}
\end{bmatrix} = \Vec{0}.
\label{eq:discrete_problem}
\end{equation}
%
Matrices $M(\Vec{s}^m)$ and $W=W^m$ are diagonal, while $D$ and $\chak{-D^T}$ resemble the discrete divergence and gradient operators respectively (note that $D$ and $D^T$ do not contain the mesh size $h$).
The nonlinear operator $T^m$ reads:
%
\begin{equation}
T^m(\Vec{\sigma}^{m}, \Vec{s}^{m}) = (\Delta t_m)^{-1}W \Vec{s}^{m} + D\diag{\Vec{\sigma}^m}U(\Vec{\sigma}^{m}) f_w(\Vec{s}^{m}),
\label{eq:discrete_transport}
\end{equation}
%
where $\diag{\Vec{\sigma}^m}$ is the diagonal matrix created from the entries of the argument vector, and $U(\Vec{\sigma})$ is the upwind operator selecting for each connection the appropriate upstream value from the input vector $f_w(\Vec{s}^{m})$, which contains the fractional flow function values evaluated in each cell.
\begin{figure}
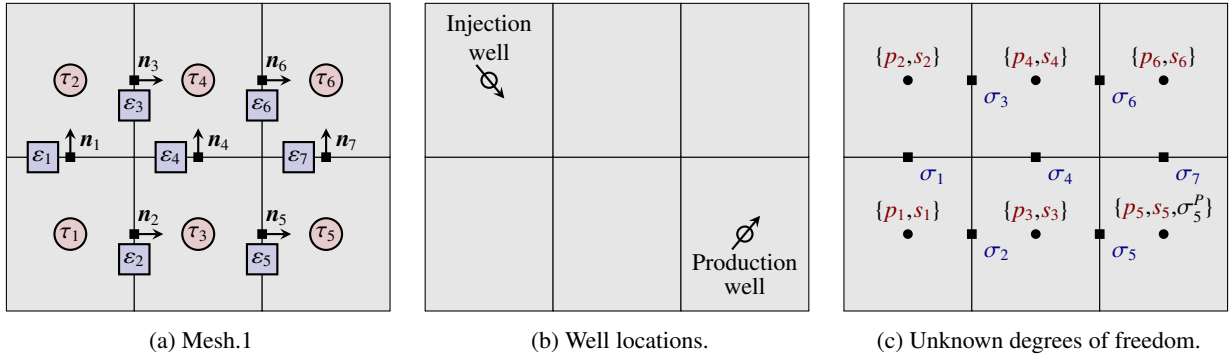

  \small
  \centering
  \begin{subfigure}[b]{0.33\textwidth}
    \centering
    \if \generateTikzFigures 1
      \include{./pics/paper_pics/app_mesh_geometry_cells}
    \else
      \includegraphics[scale=1]{main-figure\theFigureCounter.pdf}
      \stepcounter{FigureCounter}
    \fi  
    \caption{Mesh.\theFigureCounter}
    \label{fig:mesh_sketch_NEW}
  \end{subfigure}
  \hfill
  \begin{subfigure}[b]{0.33\textwidth}
    \centering
    \if \generateTikzFigures 1
      \include{./pics/paper_pics/app_mesh_geometry_edges}
    \else
      \includegraphics[scale=1]{main-figure\theFigureCounter.pdf}
      \stepcounter{FigureCounter}
    \fi
    \caption{Well locations.}
    \label{fig:mesh_well_NEW}
  \end{subfigure}
  \hfill
  \begin{subfigure}[b]{0.33\textwidth}
    \centering
    \if \generateTikzFigures 1
      \include{./pics/paper_pics/app_mesh_dofs_NEW}
    \else
      \includegraphics[scale=1]{main-figure\theFigureCounter.pdf}
      \stepcounter{FigureCounter}
    \fi
    \caption{Unknown degrees of freedom.}
    \label{fig:dof_sketch_NEW}
  \end{subfigure}
  \caption{Sketch of a well-driven flow using a mesh consisting of six cells.  The domain boundary is subject to no-flow conditions everywhere.  The location of the rate controlled injection well and the BHP controlled production well is shown in (b).}
  \label{fig:well_driven_flow}
\end{figure} 
For clarity, using a simple well-driven flow example defined in Fig. \ref{fig:well_driven_flow}, we provide additional details on the matrices and vectors appearing in the discrete residual equations \eqref{eq:discrete_problem} in Figs. \ref{fig:total_mass_residual}--\ref{fig:upwind_operator}.
\begin{figure}
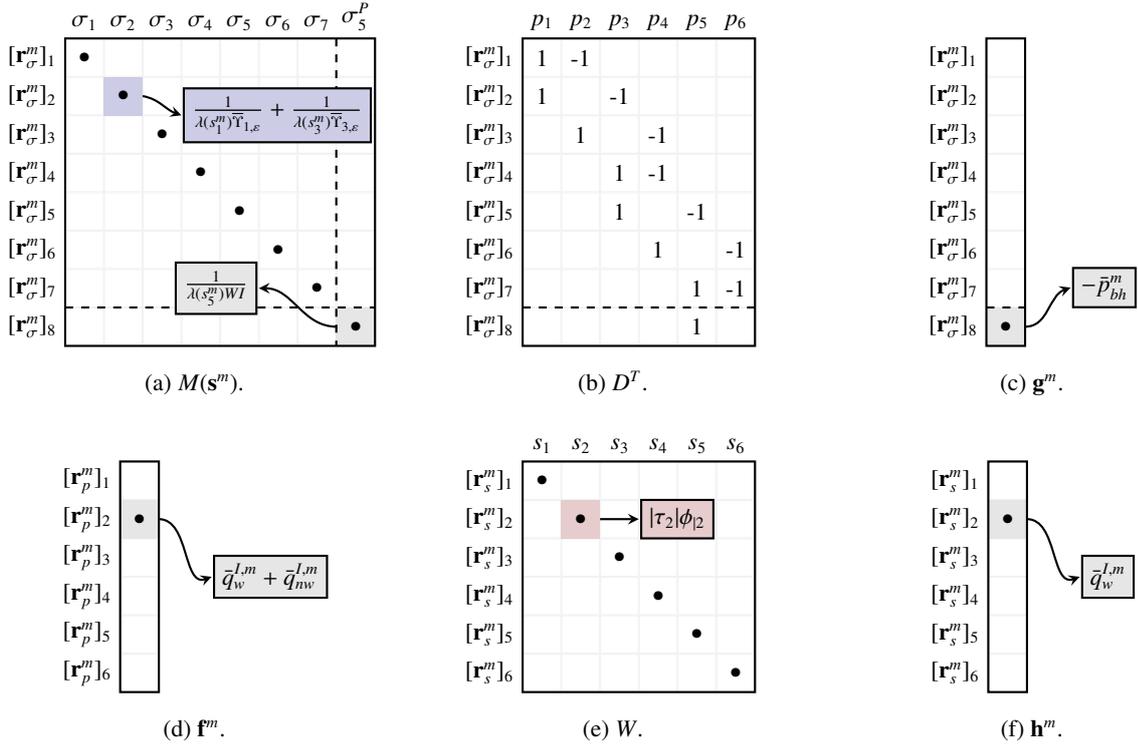

    \small
    \centering    
    \hfill
    \begin{subfigure}[b]{0.33\textwidth}
      \centering
      \if \generateTikzFigures 1
        \include{./pics/paper_pics/app_mat_M}
      \else
      \includegraphics[scale=1]{main-figure\theFigureCounter.pdf}
      \stepcounter{FigureCounter}
      \fi
      \caption{$\Mat{M}(\Vec{s}^m)$.}
    \end{subfigure}
    \hfill
    \begin{subfigure}[b]{0.33\textwidth}
      \centering
      \if \generateTikzFigures 1
        \include{./pics/paper_pics/app_mat_Dt}
      \else
      \includegraphics[scale=1]{main-figure\theFigureCounter.pdf}
      \stepcounter{FigureCounter}
      \fi
      \caption{$\Mat{D}^T$.}
    \end{subfigure} 
    \hfill
    \begin{subfigure}[b]{0.33\textwidth}
      \centering
      \if \generateTikzFigures 1
        \include{./pics/paper_pics/app_vec_g}
      \else
      \includegraphics[scale=1]{main-figure\theFigureCounter.pdf}
      \stepcounter{FigureCounter}
      \fi
      \caption{$\Vec{g}^m$.}
    \end{subfigure}    
    \hfill\null

    \hfill
    \begin{subfigure}[b]{0.33\textwidth}
      \centering
      \if \generateTikzFigures 1
        \include{./pics/paper_pics/app_vec_f}
      \else
      \includegraphics[scale=1]{main-figure\theFigureCounter.pdf}
      \stepcounter{FigureCounter}
      \fi
      \caption{$\Vec{f}^m$.}
    \end{subfigure}   
    \hfill
    \begin{subfigure}[b]{0.33\textwidth}
      \centering
      \if \generateTikzFigures 1
        \include{./pics/paper_pics/app_mat_W}
      \else
      \includegraphics[scale=1]{main-figure\theFigureCounter.pdf}
      \stepcounter{FigureCounter}
      \fi
      \caption{$\Mat{W}$.}
    \end{subfigure}  
    \hfill
    \begin{subfigure}[b]{0.33\textwidth}
      \centering
      \if \generateTikzFigures 1
        \include{./pics/paper_pics/app_vec_h}
      \else
      \includegraphics[scale=1]{main-figure\theFigureCounter.pdf}
      \stepcounter{FigureCounter}
      \fi
      \caption{$\Vec{h}^m$.}
    \end{subfigure}    
    \hfill\null
    \caption{Matrices and vectors appearing in the discrete flux and total mass residual equations \eqref{eq:discrete_problem} for the simple well-driven flow problem defined in Fig. \ref{fig:well_driven_flow}. Black dots denote a nonzero entry.}
    \label{fig:total_mass_residual}
\end{figure}
\begin{figure}
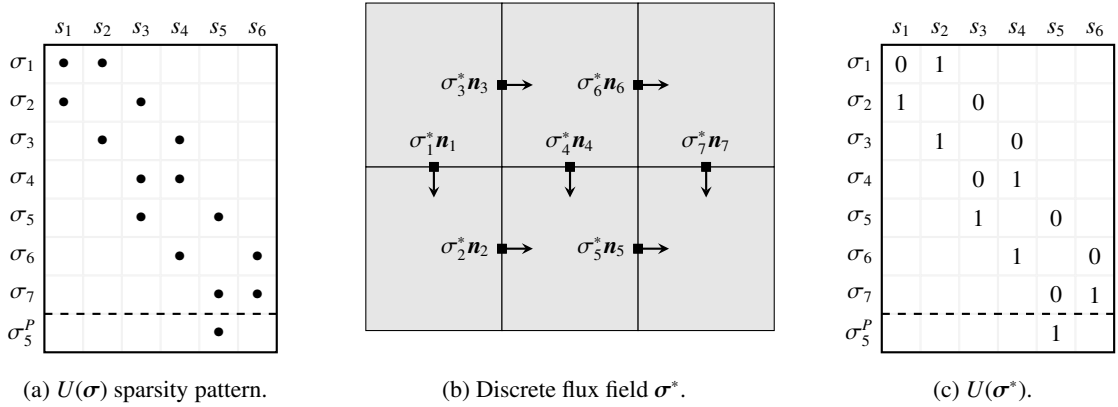

    
    \small
    \centering

    \hfill
    \begin{subfigure}[b]{0.33\textwidth}
      \centering
      \if \generateTikzFigures 1
        \include{./pics/paper_pics/app_mat_U_sparsity}
      \else
      \includegraphics[scale=1]{main-figure\theFigureCounter.pdf}
      \stepcounter{FigureCounter}
      \fi
      \caption{$\Mat{U}(\Vec{\sigma})$ sparsity pattern.}
    \end{subfigure} 
    \hfill
    \begin{subfigure}[b]{0.33\textwidth}
      \centering
      \if \generateTikzFigures 1
        \include{./pics/paper_pics/app_mat_U_flux_map}   
      \else
      \includegraphics[scale=1]{main-figure\theFigureCounter.pdf}
      \stepcounter{FigureCounter}   
      \fi
      \caption{Discrete flux field $\Vec{\sigma}^*$.}
    \end{subfigure} 
    \hfill
    \begin{subfigure}[b]{0.33\textwidth}
      \centering
      \if \generateTikzFigures 1
        \include{./pics/paper_pics/app_mat_U_example}
      \else
      \includegraphics[scale=1]{main-figure\theFigureCounter.pdf}
      \stepcounter{FigureCounter}
      \fi
      \caption{$\Mat{U}(\Vec{\sigma}^*)$.}
    \end{subfigure}    
    \hfill\null
    
    \caption{Upwind direction selection operator $U$ used in the definition of the nonlinear operator $T^m(\Vec{\sigma}^m, \Vec{s}^m)$, see Eq. \eqref{eq:discrete_transport}, for the simple well-driven flow problem defined in Fig. \ref{fig:well_driven_flow} assuming the discrete flux field $\Vec{\sigma}^*$ shown in (b). \chak{Note that $\sigma^*_1, \sigma^*_4$ and $\sigma^*_7$ are negative.}}
    \label{fig:upwind_operator}
\end{figure}

\section{Nonlinear multigrid}\label{sec:multigrid}

In this section, we propose a nonlinear multigrid solver for the discrete nonlinear system \eqref{eq:discrete_problem} that is based on the Full Approximation Scheme (FAS) \cite{brandt77,henson2003multigrid}.
We start by giving a high-level overview of FAS and its essential components.
First, we will need three intergrid transfer operators---namely, an interpolation operator $\blkInterpolate$, a restriction operator $\blkRestrict$, and a projection operator $\blkInject$.
In particular, $\blkInterpolate$ and $\blkInject$ satisfy
\begin{equation}
\blkInject\blkInterpolate = \blkMat{I}^{\ell+1},
\end{equation}
where $\blkMat{I}^{\ell+1}$ is the identity operator on the level $\ell+1$.
We use the convention that level $\ell = 0$ refers to the finest level (i.e., the original problem), and a larger value of $\ell$ means a coarser level.
Moreover, a hierarchy of nonlinear operators $\left\{ \Vec{r}^{m, \ell}(\Vec{x}^\ell) \right\}_{\ell = 0}^{L-1}$ approximating $\Vec{r}^m(\Vec{x})$ will need to be constructed.
Lastly, the approximated solution is updated at each level based on some smoothing step, denoted ``{\tt NonlinearSmoothing}".
A typical step at level $\ell$ in the full approximation scheme multigrid is stated in Algorithm~\ref{alg:fas}, where $n_s^\ell$ is the number of smoothing steps at level $\ell$.
The backtracking procedure is described in \cite[Algorithm~1]{fas-spectral-diffusion}.

The multigrid solver for \eqref{eq:discrete_problem} starts with a fine initial guess, $\blkVec{x}^{m,0} := \blkVec{x}^{m-1}$, chosen to be the converged state at the previous time step $m-1$. Then, the solver performs a sequence of nonlinear iterations denoted by the superscript $k$, as follows:
\begin{equation}
\blkVec{x}^{m,k} =  {\tt NonlinearMG}(0,\, \blkVec{x}^{m,k-1},\, \Vec{0}), \quad \forall\, k \ge 1,
\end{equation}
until a certain stopping criterion is satisfied.
In the rest of this section, the details of all multigrid cycle components will be discussed.
\begin{rmk}[Abuse of terminology]

$\blkInject$ is not a projection according to the usual definition of projections. Nevertheless, following the discussion in \cite[Remark~6]{fas-spectral-diffusion}, $\blkInject$ will be referred to as a projection with an abuse of terminology.

\end{rmk}
\begin{algorithm}[t]
\caption{Nonlinear step at level $\ell$ in the Full Approximation Scheme}
\label{alg:fas}
\begin{algorithmic}[1]
\Function{\tt NonlinearMG}{$\ell,\, \blkVec{x}^\ell,\, \blkVec{b}^\ell$}
	\If{$\ell$ is the coarsest level}
		\State $\blkVec{x}^\ell \leftarrow$ \Call{\tt NonlinearSmoothing}{$\ell,\, \blkVec{x}^\ell,\, \blkVec{b}^\ell,\, n_s^\ell$}
	\Else
		\State $\blkVec{x}^\ell \leftarrow$ \Call{\tt NonlinearSmoothing}{$\ell,\, \blkVec{x}^\ell,\, \blkVec{b}^\ell,\, n_s^\ell$} \label{alg:line:nonlinPreSmoothing}
		\State $\blkVec{x}^{\ell+1} \leftarrow \blkInject \blkVec{x}^\ell$
		\State $\blkVec{b}^{\ell+1} \leftarrow \blkVec{r}^{m, \ell+1}(\blkVec{x}^{\ell+1}) - \blkRestrict(\blkVec{r}^{m, \ell}(\blkVec{x}^\ell) - \blkVec{b}^\ell )$ 
		\State $\blkVec{y}^{\ell+1} \leftarrow$ \Call{\tt NonlinearMG}{$\ell+1,\, \blkVec{x}^{\ell+1},\, \blkVec{b}^{\ell+1}$}
		\State $\blkVec{x}^\ell \leftarrow$ \Call{\tt Backtracking}{$\blkVec{x}^\ell,\, \blkInterpolate( \blkVec{y}^{\ell+1} - \blkVec{x}^{\ell+1}),\, \theta$} 
		\State $\blkVec{x}^\ell \leftarrow$ \Call{\tt NonlinearSmoothing}{$\ell,\, \blkVec{x}^\ell,\, \blkVec{b}^\ell,\, n_s^\ell$} \label{alg:line:nonlinPostSmoothing}
  	\EndIf
	\State \Return{$\blkVec{x}^\ell$}
 \EndFunction
 \end{algorithmic}
\end{algorithm}

\subsection{Intergrid transfer operators} \label{sec:intergrid_operators}

The interpolation operator $\blkInterpolate$ and the projection operator $\blkInject$ are block-diagonal, composed of the corresponding operators for the flux, the pressure, and the saturation unknowns:
\begin{equation}
\blkInterpolate = \begin{bmatrix}
\interpolate{\sigma}  \\ & \interpolate{p} \\ & & \interpolate{s}
\end{bmatrix} 
\qquad \text{ and } \qquad 
\blkInject = \begin{bmatrix}
\inject{\sigma}  \\ & \inject{p} \\ & & \inject{s}
\end{bmatrix}.
\label{eq:prolongation_and_projection}
\end{equation}
The restriction operator $\blkRestrict$ is taken as the transpose of the interpolation operator $\blkInterpolate$, i.e., 
\[
\blkRestrict := \begin{bmatrix}
\restrict{\sigma}  \\ & \restrict{p} \\ & & \restrict{s}
\end{bmatrix} := \begin{bmatrix}
\left(\interpolate{\sigma}\right)^T  \\ & \left(\interpolate{p}\right)^T \\ & & \left(\interpolate{s}\right)^T
\end{bmatrix}.
\]

To define our interpolation operators, we first form a nested hierarchy of grids $\{\grid^\ell\}_{\ell=0}^{\mathcal{L}-1}$ by aggregating fine grid cells in $\grid^0 := \grid$.
Starting with $\ell = 0$, we consider the cell-connectivity graph of $\grid^\ell$, where each cell (respectively face) in $\grid^\ell$ is a vertex (respectively edge) of the graph.
Based on the cell-connectivity graph, contiguous aggregates of cells are formed by using a graph partitioner (e.g., METIS \cite{karypis1998fast}).
These aggregates are the ``cells" (which have irregular shapes) in the coarser-level grid $\grid^{\ell+1}$.
A coarser-level face is also naturally formed by collecting the fine faces sharing a pair of adjacent aggregates.
The set of faces on level $\ell$ is denoted by $\faces^\ell$.
This process is repeated until the coarsest grid $\grid^{\mathcal{L}-1}$ is formed.

For pressure and saturation, the coarse spaces are taken to be the space of piecewise-constant functions on the coarse grids.
The corresponding interpolation operators are defined as
\begin{equation}
\left[\interpolate{p}\right]_{ij} = \left[\interpolate{s}\right]_{ij} = \left\{ \begin{array}{ll}
1, & \text{ if aggregate $j$ in $\grid^{\ell+1}$ contains cell $i$ in $\grid^\ell$,} \\
0, & \text{ otherwise.}
\end{array}
\right.
\end{equation}
The projection operators for pressure and saturation are chosen to be
\begin{equation}
\inject{p} = \inject{s} := \left( \restrict{s} \interpolate{s} \right)^{-1} \restrict{s}.
\label{eq:simple_sat_projector}
\end{equation}
Note that, for a saturation vector $\Vec{s}^\ell$ on level $\ell$, the projection $\inject{s}\Vec{s}^\ell$ on the coarse level $\ell+1$ corresponds to taking the arithmetic average of the entries of $\Vec{s}^\ell$ in each aggregate. Also, the product $\restrict{s} \interpolate{s}$ is actually a diagonal matrix, so $\inject{s}$ is computationally cheap to obtain.
\begin{rmk}
A more physically meaningful projection would be to take a weighted average of the saturation values in each aggregate,
where the weights are the pore volume of cells.
That is,
\begin{equation}
\inject{s,W} := \left( \restrict{s} W^\ell \interpolate{s} \right)^{-1} \restrict{s} W^\ell.
\label{eq:weighted_sat_projector}
\end{equation}
In fact, we have tested both choices of the projection operators \eqref{eq:simple_sat_projector} and \eqref{eq:weighted_sat_projector}.
However, we did not see clear benefit of using one over the other in our numerical experiments.
For simplicity, we present results based on the arithmetic average.
\end{rmk}

For the coarse flux space, each basis function is associated with a coarse face $\varepsilon^{\ell+1}_i \in \faces^{\ell+1}$.
Specifically, let $\tau^{\ell+1}_K$ and $\tau^{\ell+1}_L \in \grid^{\ell+1}$ be the aggregates sharing the coarse face $\varepsilon^{\ell+1}_i = \varepsilon^{\ell+1}_{K, L} := \partial\tau^{\ell+1}_K \cap \partial\tau^{\ell+1}_L$.
The basis function associated with $\varepsilon^{\ell+1}_i$ is obtained by first solving a local boundary value problem in $\tau^{\ell+1}_K \cup \tau^{\ell+1}_L$ discretized on level $\ell$:
\begin{equation}
\begin{split}
         \mathbb{K}^{-1} \cdot \widetilde{\tensorOne{\phi}}_i^{\ell+1} + \nabla p & = 0 \\
	\nabla \cdot \widetilde{\tensorOne{\phi}}_i^{\ell+1} & = q_{\tau^{\ell+1}_K \cup \tau^{\ell+1}_L}
\end{split}
\label{eq:local_flux_problem}
\end{equation}
with no-flow boundary condition $\widetilde{\tensorOne{\phi}}_i^{\ell+1} \cdot \tensorOne{n} = 0$ on $\partial(\tau^{\ell+1}_K \cup \tau^{\ell+1}_L)$,
where
\begin{equation}
q_{\tau^{\ell+1}_K \cup \tau^{\ell+1}_L} = \left\{ 
\begin{array}{ll}
1 / |\tau^{\ell+1}_K|, & \text{ in } \tau^{\ell+1}_K, \\
-1 / |\tau^{\ell+1}_L|, & \text{ in } \tau^{\ell+1}_L.
\end{array}
\right.
\label{eq:flux_normalization}
\end{equation}
Note that $\widetilde{\tensorOne{\phi}}_i^{\ell+1}$ was introduced in mixed multiscale finite element methods; see, for example, \cite{aarnes06, arbogast2012}.
The final basis $\tensorOne{\phi}_i^{\ell+1}$ is obtained by normalizing $\widetilde{\tensorOne{\phi}}_i^{\ell+1}$ so that the total normal flux of $\tensorOne{\phi}_i^{\ell+1}$ on $\varepsilon^{\ell+1}_i$ is 1. More precisely,
\begin{equation}
\tensorOne{\phi}_i^{\ell+1} := \left( \int_{\varepsilon^{\ell+1}_i}  \widetilde{\tensorOne{\phi}}_i^{\ell+1} \cdot \tensorOne{n} \right)^{-1} \widetilde{\Vec{\phi}}_i^{\ell+1}.
\end{equation}
The interpolation operator $\interpolate{\sigma}$ is formed by collecting the coefficient vectors of the local flux solutions (corresponding to $\tensorOne{\phi}_i^{\ell+1} $) as its column vectors. 
The projection $\inject{\sigma}$ is also defined locally on each coarse face. 
For a given finer level flux $\sigma^\ell$, the $i$-th entry of $\inject{\sigma} \sigma^\ell$ is the total normal flux of $\sigma^\ell$ on the coarse face $\varepsilon^{\ell+1}_i$.
In the next subsection, we will introduce discrete operators $M^\ell(\Vec{s})$ and $D^\ell$, which are the counterparts of $M(\Vec{s})$ and $D$ on level $\ell$.
In practice, the discrete problems of \eqref{eq:local_flux_problem} and \eqref{eq:flux_normalization},
as well as $\inject{\sigma}$, are constructed using submatrices of $M^\ell(\Vec{0})$ and $D^\ell$; cf. Appendices A and B of \cite{ml-spectral-coarsening}. 
We remark here that the coarsening in the current paper is the lowest order version of the spectral coarsening method \cite{ml-spectral-coarsening, fas-spectral-diffusion}.

\begin{rmk}
\label{remark:wells}
In this work, we do not consider the impact of wells during the coarse space construction.  For simplicity, cells with well perforations are kept as separate aggregates and carried unmodified to coarser levels.  A specific treatment accounting for well interactions is the subject of ongoing work.  Given that wells drive the flow dynamics, a more sophisticated coarse space could offer performance advantages.
\end{rmk}

\subsection{The nonlinear problem on each level} \label{sec:coarse_operators}

On the fine level $\ell = 0$, let $\Vec{r}^{m, 0}(\Vec{x}^0) := \Vec{r}^m(\Vec{x}^0)$, $M^0(\Vec{s}^0) := M(\Vec{s}^0)$,
$D^0 := D$, $T^{m,0}(\Vec{\sigma}^0, \Vec{s}^0) := T^m(\Vec{\sigma}^0, \Vec{s}^0)$, $\Vec{g}^{m, 0} := \Vec{g}^{m}$, 
$\Vec{f}^{m, 0} := \Vec{f}^{m}$, and $\Vec{h}^{m, 0} := (\Delta t_m)^{-1}W\Vec{s}^{m-1} + \Vec{h}^{m}$. 
With the interpolation operator $\blkInterpolate$ and restriction operator $\blkRestrict$, the nonlinear operators on coarse levels are defined recursively as
\begin{equation}
\Vec{r}^{m, \ell+1}(\Vec{x}^{\ell+1}) := \blkRestrict\Vec{r}^{m, \ell}\left( \blkInterpolate\Vec{x}^{\ell+1} \right)
= 
\begin{bmatrix}
M^{\ell+1}\left(\Vec{s}^{\ell+1} \right)\Vec{\sigma}^{\ell+1} \chak{-(D^{\ell+1})^T} \Vec{p}^{\ell+1} - \Vec{g}^{m,\ell+1}\\
D^{\ell+1} \Vec{\sigma}^{\ell+1} - \Vec{f}^{m,\ell+1}\\
T^{m, \ell+1}(\Vec{\sigma}^{\ell+1}, \Vec{s}^{\ell+1}) - \Vec{h}^{m,\ell+1}
\end{bmatrix}
=: 
\begin{bmatrix}
\Vec{r}_\sigma^{m, \ell+1}(\Vec{\sigma}^{\ell+1}, \Vec{p}^{\ell+1}, \Vec{s}^{\ell+1}) \\
\Vec{r}_p^{m, \ell+1}(\Vec{\sigma}^{\ell+1}) \\
\Vec{r}_s^{m, \ell+1}(\Vec{\sigma}^{\ell+1}, \Vec{s}^{\ell+1})  
\end{bmatrix}
\quad
\label{eq:level_problem}
\end{equation}
where
\begin{subequations}
\begin{align}
M^{\ell+1}\left( \Vec{s}^{\ell+1} \right) & := \restrict{\sigma} M^\ell \left(  \interpolate{s}\Vec{s}^{\ell+1} \right) \interpolate{\sigma},
& 
\Vec{g}^{m,\ell+1} &:=  \restrict{\sigma} \Vec{g}^{m,\ell}, \nonumber \\
D^{\ell+1} & :=  \restrict{p} D^\ell \interpolate{\sigma},
&
\Vec{f}^{m,\ell+1} &:=  \restrict{p} \Vec{f}^{m,\ell},  \nonumber \\
T^{m, \ell+1}(\Vec{\sigma}^{\ell+1}, \Vec{s}^{\ell+1}) & :=  \restrict{s} T^{m, \ell}\left( \interpolate{\sigma} \Vec{\sigma}^{\ell+1},  \interpolate{s} \Vec{s}^{\ell+1} \right), 
&
\Vec{h}^{m,\ell+1} &:=  \restrict{s} \Vec{h}^{m,\ell}. \nonumber
\end{align}
\end{subequations}
Note that \eqref{eq:level_problem} is a conceptual definition of the coarse operators. 
In practice, due to scalability concerns, we do not want the evaluation of $\Vec{r}^{m, \ell+1}(\Vec{x}^{\ell+1})$ during the multigrid cycle to involve computations on the finer level $\ell$.
To achieve this, we have to construct and store some coarse operators during the setup phase of the multigrid solver.
For $D^{\ell+1}, \Vec{g}^{m,\ell+1}, \Vec{f}^{m,\ell+1}$, and $\Vec{h}^{m,\ell+1}$, the construction is straightforward.
The main issue is in the evaluation of the nonlinear components $M^{\ell+1}\left( \Vec{s}^{\ell+1} \right)$ and $T^{m, \ell+1}(\Vec{\sigma}^{\ell+1}, \Vec{s}^{\ell+1})$.
In the rest of this section, we will describe how it can be done efficiently.

\subsubsection{Evaluation of $M^{\ell+1}\left( \Vec{s}^{\ell+1} \right)$}

For $M^{\ell+1}\left( \Vec{s}^{\ell+1} \right)$, we follow the procedure proposed in \cite[Section 3.2]{fas-spectral-diffusion}.
To begin with, some local matrices $\widehat{M}^{\ell}_{\tau^{\ell}_K}$ that are independent of $\lambda(\Vec{s}^\ell)$,
each of which is associated with a cell $\tau^{\ell}_K\in \grid^{\ell}$, are precomputed in the setup phase.
On the fine level, $\widehat{M}^{0}_{\tau^{0}_K}$ is a diagonal matrix whose entries are all the half transmissibility associated with the faces on the boundary of $\tau^{0}_K$:
\begin{equation}
\widehat{M}^{0}_{\tau^{0}_K} = \begin{bmatrix}
\ddots & \\
& \overline{\Upsilon}_{K,\varepsilon^0} \\
& & \ddots
\end{bmatrix}_{\varepsilon^0\in \partial \tau^{0}_K}
\end{equation}
where the half transmissibility $\overline{\Upsilon}_{K,\varepsilon^0}$ is defined in \eqref{eq:half-transmissibility}.
On coarse levels, $\widehat{M}^{\ell+1}_{\tau^{\ell+1}_K}$ is obtained by first assembling the finer-level local matrices $\widehat{M}^{\ell}_{\tau^{\ell}_L}$ associated with the fine cells $\tau^{\ell}_L$ covering the coarse cell $\tau^{\ell+1}_K$,
followed by a local variational coarsening using $\interpolate{\sigma}$.

Now let $\Vec{s}^{\ell+1}$ be a coarse saturation on level $\ell+1$.
By our construction, the $K$-th entry of $\Vec{s}^{\ell+1}$, $s^{\ell+1}_K$, represents the average saturation value in the coarse cell $\tau^{\ell+1}_K\in \grid^{\ell+1}$.
During the solving phase, $M^{\ell+1}\left( \Vec{s}^{\ell+1} \right)$ is obtained by assembling the precomputed local matrix $\widehat{M}^{\ell+1}_{\tau^{\ell+1}_K}$ scaled by $\lambda(s^{\ell+1}_K)^{-1}$ using a local-to-global map.
Notice that the assembling of $M^{\ell+1}\left( \Vec{s}^{\ell+1} \right)$ in this approach involves quantities on the coarse level $\ell+1$ only.

\subsubsection{Evaluation of $\Mat{T}^{m, \ell}(\Vec{\sigma}^{\ell+1}, \Vec{s}^{\ell+1})$}

Let $W^0 := W$ and $U^0(\Vec{\sigma}^0) := U(\Vec{\sigma}^0)$. Note that $\Mat{T}^{m, 0}$ has two parts:
\begin{equation}
\Mat{T}^{m, 0}(\Vec{\sigma}^0, \Vec{s}^0) :=  \Mat{T}^{m, 0}_1(\Vec{s}^0) + \Mat{T}^{m, 0}_2(\Vec{\sigma}^0, \Vec{s}^0) := (\Delta t_m)^{-1}W^0 \Vec{s}^0 + D^0\diag{\Vec{\sigma}^0}U^0(\Vec{\sigma}^0) f_w(\Vec{s}^0).
\label{eq:T_split}
\end{equation}
The coarsening of $\Mat{T}^{m, 0}_1(\Vec{s}^0)$ is straightforward as it is a linear operator:
\[
\Mat{T}^{m, \ell+1}_1(\Vec{s}^{\ell+1}) := \restrict{s}\Mat{T}^{m, \ell}_1(\interpolate{s}\Vec{s}^{\ell+1}) := (\Delta t_m)^{-1} W^{\ell+1} \Vec{s}^{\ell+1}
\]
where $W^{\ell+1}$ is the mass matrix for the saturation space on the coarse level $\ell+1$ defined as:
\begin{equation}
W^{\ell+1} :=  \restrict{s} W^\ell \interpolate{s}, \; \forall \ell \ge0.
\end{equation}

The coarsening of $\Mat{T}^{m, 0}_2$ is more challenging as it involves a nonlinear function and upwind fluxes.
Nevertheless, we can actually compute $\Mat{T}^{m, \ell+1}_2(\Vec{\sigma}^{\ell+1}, \Vec{s}^{\ell+1})$ without visiting the finer levels during the multigrid solving phase.
To see this, consider the two-level case.
In this case, we use some lighter notation to simplify the presentation.
Specifically, the superscript and subscript in intergrid operators are dropped.
For coefficient vectors, operators, and geometrical entities on the fine level, we drop the superscript 0.
On the other hand, for coefficient vectors, operators, and geometrical entities on the coarse level, the superscript $1$ is replaced by $c$.
For example, $\interpolateOne{\sigma} :=  \big(\Mat{P}_\sigma \big)_1^{0}$, $\Vec{\sigma} = \Vec{\sigma}^0$, $\Vec{s}^c := \Vec{s}^1$, $\Mat{T}^{m}_2 := \Mat{T}^{m, 0}_2$, $\Mat{D}^c := \Mat{D}^1$, and $\grid^c := \grid^1$.
The following lemma characterizes the coarse upwind fluxes when saturation is piecewise constant on the coarse level.
\begin{lemma}
\label{lemma:coarse_upwind_flux_characterization}
Let $\Vec{s}^c = (s^c_K)_{\tau^c_K \in \grid^c}$ and $\Vec{\sigma} = (\sigma_j)_{\varepsilon_j \in \faces}$ be the coefficient vectors of some coarse saturation solution and fine flux solution respectively.
Then, the $K$-th entry of the variationally defined quantity  $\restrictOne{s}\Mat{T}^{m}_2\left(\Vec{\sigma}, \interpolateOne{s}\Vec{s}^c \right)$ reads as 
\begin{equation}
\left( \restrictOne{s}\Mat{T}^{m}_2\left(\Vec{\sigma}, \interpolateOne{s}\Vec{s}^c \right) \right)_K
= \sum_{\varepsilon^c_i = \varepsilon^c_{K, L} \subseteq \partial\tau^c_K} \left( f_w(s^c_K) \sum_{\varepsilon_j \subseteq \varepsilon^c_i} \max\left( \sigma_j, 0 \right)  + f_w(s^c_{L})  \sum_{\varepsilon_j \subseteq \varepsilon^c_i} \min\left( \sigma_j, 0 \right)  \right)
\label{eq:coarse_upwind_characterization}
\end{equation}
where $\varepsilon^c_i = \varepsilon^c_{K, L}$ is the coarse face shared by the coarse cells $\tau^c_K$ and $\tau^c_{L}$.
\end{lemma}
\begin{proof}
First of all, since our coarse saturation space consists of piecewise constant functions, the order of evaluation of nonlinear function and interpolation can be switched:
\begin{equation}
f_w \left( \interpolateOne{s} \, \Vec{s}^c \right) = \interpolateOne{s}  f_w \left( \Vec{s}^c \right).
\label{eq:nonlinear_Ps_commute}
\end{equation}
Let $\Vec{\sigma}_{f^{upw}} := \diag{ \Vec{\sigma} } U \left(  \Vec{\sigma} \right) \interpolateOne{s} f_w \left( \Vec{s}^c \right)$.
By the definition of $\Mat{T}^{m}_2$ given in \eqref{eq:T_split}, and \eqref{eq:nonlinear_Ps_commute},
\begin{equation}
\restrictOne{s}\Mat{T}^{m}_2\left(\Vec{\sigma}, \interpolateOne{s}\Vec{s}^c \right)
= \restrictOne{s} D \, \diag{ \Vec{\sigma} } U \left(  \Vec{\sigma} \right) \interpolateOne{s} f_w \left( \Vec{s}^c \right)
= \restrictOne{s}D\Vec{\sigma}_{f^{upw}}.
\label{eq:coarse_upwind_step0}
\end{equation}
Note that $\Vec{\sigma}_{f^{upw}}$ is just another fine-level flux.
So $D\Vec{\sigma}_{f^{upw}}$ computes the sum of outward normal flux (of the flux represented by $\Vec{\sigma}_{f^{upw}}$) of each fine-level cell.
By the definition of $\restrictOne{s}$, for each coarse cell $\tau^c_K$, 
$\restrictOne{s}D\Vec{\sigma}_{f^{upw}}$ sums up the outward normal fluxes of the fine-level cells $\tau_L$ in the coarse cell $\tau^c_K$.
Since $\Vec{\sigma}_{f^{upw}}$ is conservative, for two cells sharing a common face, outward normal fluxes of the cells on the common face cancel each other. 
This means that the $K$-th entry of $\restrictOne{s}D\Vec{\sigma}_{f^{upw}}$ is the sum of outward normal flux of the coarse-level cell $\tau^c_K$.

Now for the fine faces $\varepsilon_j$ on $\partial\tau^c_K$, we can group them based on the coarse face they belong to, and look at one coarse face $\varepsilon^c_i$ at a time.
Let $\varepsilon^c_i = \varepsilon^c_{K, L}$ be the coarse face shared by the coarse cells $\tau^c_K$ and $\tau^c_{L}$.
Notice that $\Vec{\sigma}_{f^{upw}}$ is actually some upwind flux on the fine level.
Since $\interpolateOne{s} f_w \left( \Vec{s}^c \right)$ is piecewise constant on coarse level,
for the fine faces belonging to a coarse face, there are only two possible upwinded quantities: $f_w(s^c_K)$ or $f_w(s^c_L)$, depending on the direction of $\Vec{\sigma}$ on the fine faces.
Hence, the $K$-th entry of $\restrictOne{s}D\Vec{\sigma}_{f^{upw}}$ is
\begin{equation}
\begin{split}
\left( \restrictOne{s}D\Vec{\sigma}_{f^{upw}} \right)_K
= & \sum_{\varepsilon^c_i = \varepsilon^c_{K, L} \subseteq \partial\tau^c_K}  \left( \sum_{\varepsilon_j \subseteq \varepsilon^c_i} \left( \max\left( \sigma_j, 0 \right) f_w(s^c_K) + \min\left( \sigma_j, 0 \right) f_w(s^c_{L})  \right) \right) \\
= & \sum_{\varepsilon^c_i = \varepsilon^c_{K, L} \subseteq \partial\tau^c_K} \left( f_w(s^c_K) \sum_{\varepsilon_j \subseteq \varepsilon^c_i} \max\left( \sigma_j, 0 \right)  + f_w(s^c_{L})  \sum_{\varepsilon_j \subseteq \varepsilon^c_i} \min\left( \sigma_j, 0 \right)  \right).
\end{split}
\label{eq:coarse_upwind_step1}
\end{equation}
Combining \eqref{eq:coarse_upwind_step0} and \eqref{eq:coarse_upwind_step1}, we get \eqref{eq:coarse_upwind_characterization}.

\end{proof}

Note that Lemma~\ref{lemma:coarse_upwind_flux_characterization} holds for any $\Vec{\sigma}$ on the fine level.
In particular, $\Vec{\sigma}$ can be taken to be the interpolation of $\Vec{\sigma}^c$: $\Vec{\sigma} = \interpolateOne{\sigma}\Vec{\sigma}^c$.
Therefore, by Lemma~\ref{lemma:coarse_upwind_flux_characterization}, we know that $T^{m, c}_2(\Vec{\sigma}^c, \Vec{s}^c) = \restrictOne{s}\Mat{T}^{m}_2\left(\interpolateOne{\sigma}\Vec{\sigma}^c, \interpolateOne{s}\Vec{s}^c \right) $ only depends on the values of $\interpolateOne{\sigma}\Vec{\sigma}^c$ on (the fine faces in) coarse faces.
Using this observation, and a carefully defined coarse version of $U$, we will see in the next proposition that $\Mat{T}^{m, c}_2$ actually has a similar structure as $\Mat{T}^{m}_2$.

\begin{figure}
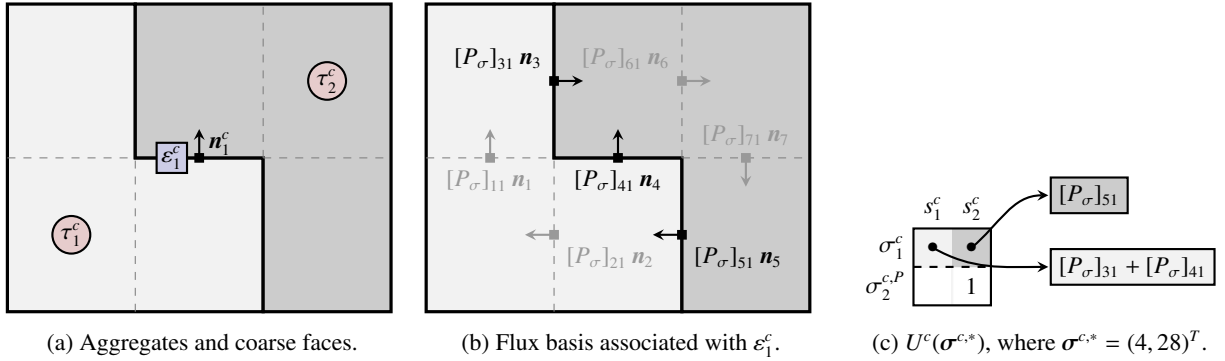

  \small
  \centering
  \begin{subfigure}[b]{0.33\textwidth}
    \centering
    \if \generateTikzFigures 1
      \include{./pics/paper_pics/app_mesh_geometry_aggregates2}
    \else
      \includegraphics[scale=1]{main-figure\theFigureCounter.pdf}
      \stepcounter{FigureCounter}
    \fi
    \caption{Aggregates and coarse faces.}
    \label{fig:mesh_sketch_aggregates_a}
  \end{subfigure}
  \hfill
  \begin{subfigure}[b]{0.33\textwidth}
    \centering
    \if \generateTikzFigures 1
      \include{./pics/paper_pics/app_flux_basis2}
    \else    
      \includegraphics[scale=1]{main-figure\theFigureCounter.pdf}
      \stepcounter{FigureCounter}
    \fi
    \caption{Flux basis associated with $\varepsilon^c_1$.}
    \label{fig:mesh_sketch_aggregates_b}
  \end{subfigure}
  \hfill
  \begin{subfigure}[b]{0.33\textwidth}
    \centering
    \if \generateTikzFigures 1
      \include{./pics/paper_pics/app_mat_Uc_example2}
    \else
      \includegraphics[scale=1]{main-figure\theFigureCounter.pdf}
      \stepcounter{FigureCounter}
    \fi
    \caption{$U^c(\Vec{\sigma}^{c,*})$, where $\Vec{\sigma}^{c,*} = (4, 28)^T$.}
  \end{subfigure}
  \caption{Coarse upwind direction selection operator $U^c$ defined in \eqref{eq:coarse_upwind_selector} for a coarsened version of the simple well-driven flow problem defined in Fig.~\ref{fig:well_driven_flow} assuming the discrete flux field of the flux basis associated with $\varepsilon^c_1$ shown in (b).}
  \label{fig:coarse_upwind_selector}
\end{figure}

\begin{proposition}
\label{prop:coarse_transport_structure}
Let $\Vec{s}^c = (s^c_K)_{\tau^c_K \in \grid^c}$ and $\Vec{\sigma}^c = (\sigma^c_i)_{\varepsilon^c_i \in \faces^c}$ be the coefficient vectors of some coarse saturation and flux solution respectively. 
Moreover, using $\Vec{\sigma}^c$, define $U^c(\Vec{\sigma}^c)$ to be the matrix such that:
\begin{align}
[U^c(\Vec{\sigma}^c)]_{ij}
& =
\begin{dcases}
  \; \delta_{jK}\sum_{\varepsilon_k \subseteq \varepsilon^c_i} \max\left(  \left[ \interpolateOne{\sigma} \right]_{ki}, 0 \right) +  \delta_{jL} \sum_{\varepsilon_k \subseteq \varepsilon^c_i} \min\left(  \left[ \interpolateOne{\sigma} \right]_{ki}, 0 \right),  & \text{if } \sigma^c_i > 0,   \\
  \; \delta_{jK}\sum_{\varepsilon_k \subseteq \varepsilon^c_i} \min\left(  \left[ \interpolateOne{\sigma} \right]_{ki}, 0 \right) +  \delta_{jL} \sum_{\varepsilon_k \subseteq \varepsilon^c_i} \max\left(  \left[ \interpolateOne{\sigma} \right]_{ki}, 0 \right),   & \text{if } \sigma^c_i \le 0,
\end{dcases}
\label{eq:coarse_upwind_selector}
\end{align}
where $\varepsilon^c_i=\varepsilon^c_{K,L}$, and $\interpolateOne{\sigma}$ is the interpolation matrix for the flux spaces. Also, $\delta_{jK}$  and $\delta_{jL}$ are the standard Kronecker symbols.
Then,
\begin{equation}
T^{m, c}_2(\Vec{\sigma}^c, \Vec{s}^c) = D^c \diag{\Vec{\sigma}^c}U^c(\Vec{\sigma}^c) f_w(\Vec{s}^c).
\end{equation}
\end{proposition}
\begin{proof}
Consider the interpolation of $\Vec{\sigma}^c$ on the fine level: $\Vec{\sigma} = \interpolateOne{\sigma}\Vec{\sigma}^c$.
By our construction, there is only one coarse flux basis vector per coarse face.
Thus, the values of $\sigma_j$ on fine faces belonging to the coarse face $\varepsilon^c_i$ is solely depending on the value of $\sigma^c_i$.
More precisely, for $i$ such that $\varepsilon_j \subseteq \varepsilon^c_i$, we have
\begin{equation}
\sigma_j = \left[ \interpolateOne{\sigma} \right]_{ji} \sigma^c_i.
\label{eq:fine_flux_on_coarse_face}
\end{equation}
In what follows in equation \eqref{eq:coarse_upwind_step3}, we use the elementary identity valid for any real numbers $p,q,a,$ and $b$
\begin{equation*}
p\max(ab,0) +q\min(ab,0)=
p\max(a,0)\max(b,0) + q\max(a,0)\min(b,0)+p\min(a,0)\min(b,0)+q\min(a,0)\max(b,0).
\end{equation*}
By Lemma~\ref{lemma:coarse_upwind_flux_characterization}, \eqref{eq:fine_flux_on_coarse_face}, and \eqref{eq:coarse_upwind_selector}, we get
\begin{equation}
\begin{split}
\left( T^{m, c}_2(\Vec{\sigma}^c, \Vec{s}^c) \right)_K
= & \; \restrictOne{s}\Mat{T}^{m}_2\left(\Vec{\sigma}, \interpolateOne{s}\Vec{s}^c \right) \\ 
= & \sum_{\varepsilon^c_i = \varepsilon^c_{K, L} \subseteq \partial\tau^c_K} \left( f_w(s^c_K) \sum_{\varepsilon_j \subseteq \varepsilon^c_i} \max\left( \left[ \interpolateOne{\sigma} \right]_{ji} \sigma^c_i, 0 \right)  + f_w(s^c_{L})  \sum_{\varepsilon_j \subseteq \varepsilon^c_i} \min\left(  \left[ \interpolateOne{\sigma} \right]_{ji} \sigma^c_i, 0 \right)  \right) \\
= & \sum_{\varepsilon^c_i = \varepsilon^c_{K, L} \subseteq \partial\tau^c_K} \max\left( \sigma^c_i, 0 \right) \left( f_w(s^c_K) \sum_{\varepsilon_j \subseteq \varepsilon^c_i} \max\left(  \left[ \interpolateOne{\sigma} \right]_{ji}, 0 \right)  + f_w(s^c_{L})  \sum_{\varepsilon_j \subseteq \varepsilon^c_i} \min\left(  \left[ \interpolateOne{\sigma} \right]_{ji}, 0 \right)  \right) \\
& + \min\left( \sigma^c_i, 0 \right) \left( f_w(s^c_K) \sum_{\varepsilon_j \subseteq \varepsilon^c_i} \min\left(  \left[ \interpolateOne{\sigma} \right]_{ji}, 0 \right)  + f_w(s^c_{L})  \sum_{\varepsilon_j \subseteq \varepsilon^c_i} \max\left(  \left[ \interpolateOne{\sigma} \right]_{ji}, 0 \right)  \right) \\
= & \sum_{\varepsilon^c_i= \varepsilon^c_{K, L}  \subseteq \partial\tau^c_K} \sigma^c_i \left( [U^c]_{iK} f_w(s^c_K) + [U^c]_{iL} f_w(s^c_L) \right) \\
= & \Big( D^c \diag{\Vec{\sigma}^c}U^c(\Vec{\sigma}^c) f_w(\Vec{s}^c) \Big)_K.
\end{split}
\label{eq:coarse_upwind_step3}
\end{equation}

\end{proof}
An example $U^c(\Vec{\sigma}^c)$ for a coarsened problem of the simple well-driven flow problem defined in Fig.~\ref{fig:well_driven_flow} is illustrated in Fig.~\ref{fig:coarse_upwind_selector}.

Lastly, we go back to the general case where we can have multiple coarse levels.
Note that because the space of saturation in a coarser level $\ell > 1$ is still piecewise constant,
and the coarse flux space also has only one basis per coarse face,
Lemma~\ref{lemma:coarse_upwind_flux_characterization} and Proposition~\ref{prop:coarse_transport_structure} can be extended to any coarse level $\ell \ge 1$.
To this end, define a special interpolation matrix that goes all the way from level $\ell+1$ to level 0:
\begin{equation}
\big(\Mat{P}_\sigma \big)_{\ell+1}^0 := \big(\Mat{P}_\sigma \big)_1^0 \big(\Mat{P}_\sigma \big)_2^1\cdots \interpolate{\sigma}.
\end{equation}
Then, given a coarse flux solution $\Vec{\sigma}^{\ell+1} = (\sigma^{\ell+1}_i)_{\varepsilon^{\ell+1}_i \in \faces^{\ell+1}}$ on level $\ell+1$, where $\varepsilon^{\ell+1}_i=\varepsilon^{\ell+1}_{K,L}$, define $U^{\ell+1}(\Vec{\sigma}^{\ell+1})$ to be
\begin{align}
[U^{\ell+1}(\Vec{\sigma}^{\ell+1})]_{ij}
& =
\begin{dcases}
  \; \delta_{jK}\sum_{\varepsilon_k \subseteq \varepsilon^{\ell+1}_i} \max\left(  \left[ \big(\Mat{P}_\sigma \big)_{\ell+1}^0 \right]_{ki}, 0 \right) +  \delta_{jL} \sum_{\varepsilon_k \subseteq \varepsilon^{\ell+1}_i} \min\left(  \left[ \big(\Mat{P}_\sigma \big)_{\ell+1}^0 \right]_{ki}, 0 \right),  & \text{if } \sigma^{\ell+1}_i > 0,   \\
  \; \delta_{jK}\sum_{\varepsilon_k \subseteq \varepsilon^{\ell+1}_i} \min\left(  \left[ \big(\Mat{P}_\sigma \big)_{\ell+1}^0 \right]_{ki}, 0 \right) +  \delta_{jL} \sum_{\varepsilon_k \subseteq \varepsilon^{\ell+1}_i} \max\left(  \left[ \big(\Mat{P}_\sigma \big)_{\ell+1}^0 \right]_{ki}, 0 \right),   & \text{if } \sigma^{\ell+1}_i \le 0.
\end{dcases}
\label{eq:general_upwind_selector}
\end{align}
Using a similar derivation as in Lemma~\ref{lemma:coarse_upwind_flux_characterization} and Proposition~\ref{prop:coarse_transport_structure}, we can conclude that
\begin{equation}
T^{m, \ell+1}(\Vec{\sigma}^{\ell+1}, \Vec{s}^{\ell+1}) = (\Delta t_m)^{-1} W^{\ell+1} \Vec{s}^{\ell+1} + D^{\ell+1} \diag{\Vec{\sigma}^{\ell+1}}U^{\ell+1}(\Vec{\sigma}^{\ell+1}) f_w(\Vec{s}^{\ell+1}).
\label{eq:T_structure}
\end{equation}

\begin{rmk}

Recall that by construction \eqref{eq:flux_normalization} the total normal flux of the coarse basis $\tensorOne{\phi}_i^{\ell+1}$ on $\varepsilon^{\ell+1}_i$ is 1.
This means that
\[
 \sum_{\varepsilon_k \subseteq \varepsilon^{\ell+1}_i} \max\left(  \left[\big(\Mat{P}_\sigma \big)_{\ell+1}^0 \right]_{ki}, 0 \right) + \sum_{\varepsilon_k \subseteq \varepsilon^{\ell+1}_i}  \min\left(  \left[ \big(\Mat{P}_\sigma \big)_{\ell+1}^0 \right]_{ki}, 0 \right)
 = \sum_{\varepsilon_k \subseteq \varepsilon^{\ell+1}_i} \left[\big(\Mat{P}_\sigma \big)_{\ell+1}^0 \right]_{ki} = 1.
\]
Hence, by letting $\Mat{P}_{\sigma, i} ^{-ve, \ell+1}$ be the sum of negative fine fluxes on $\varepsilon^{\ell+1}_i$:
\[
\Mat{P}_{\sigma, i}^{-ve, \ell+1} =  \sum_{\varepsilon_k \subseteq \varepsilon^{\ell+1}_i} \min\left(  \left[\big(\Mat{P}_\sigma \big)_{\ell+1}^0 \right]_{ki}, 0 \right) \le 0, 
\]
we can rewrite the coarse operator $U^{\ell+1}(\Vec{\sigma}^{\ell+1})$ defined in \eqref{eq:general_upwind_selector} to be
\begin{align}
[U^{\ell+1}(\Vec{\sigma}^{\ell+1})]_{ij}
& =
\begin{dcases}
  \; \delta_{jK} (1- \Mat{P}_{\sigma, i}^{-ve, \ell+1}) +  \delta_{jL} \Mat{P}_{\sigma, i}^{-ve, \ell+1} ,  & \text{if } \sigma^{\ell+1}_i > 0,   \\
  \; \delta_{jK}  \Mat{P}_{\sigma, i}^{-ve, \ell+1} +  \delta_{jL} (1 - \Mat{P}_{\sigma, i}^{-ve, \ell+1}),   & \text{if } \sigma^{\ell+1}_i \le 0.
\end{dcases}
\label{eq:simplified_general_upwind_selector}
\end{align}
It is easy to see from \eqref{eq:simplified_general_upwind_selector} that the definition of $U^{\ell+1}$ coincides with the usual upwind direction selection operator $U$ defined in Section~\ref{sec:FV_tpfa} if $\Mat{P}_{\sigma, i}^{-ve, \ell+1} = 0$.
Therefore, $U^{\ell+1}$ can be seen as a generalization of the usual upwind direction selection operator.
Moreover, only $\Mat{P}_{\sigma, i}^{-ve, \ell+1}$ needs to be stored for each coarse face $\varepsilon^{\ell+1}_i$ in the setup phase of the multigrid in order to assemble $U^{\ell+1}(\Vec{\sigma}^{\ell+1})$ later during the solving phase.

\end{rmk}


\subsection{Nonlinear smoothing}\label{sec:newton}

The nonlinear smoothing at each level $\ell$ is used to approximate the solution to the discrete problem $\Vec{r}^{m, \ell}(\Vec{x}^\ell) = \Vec{b}^\ell$. 
In our multigrid approach, the approximation is made using Newton's method.
When {\tt NonlinearSmoothing}($\ell,\, \blkVec{x}^\ell,\, \blkVec{b}^\ell,\, n_s^\ell$) in Algorithm~\ref{alg:fas} is called, we start the Newton iterations with $\Vec{x}^{\ell, 0} := \Vec{x}^{\ell}$, and unless some stopping criteria are reached, the Newton iterations continue as 
\begin{equation}
\Vec{x}^{\ell, k} := \Vec{x}^{\ell, k-1} - \left(\partial \Vec{r}^{m,\ell}( \Vec{x}^{\ell, k-1} ) \right)^{-1} \left( \Vec{r}^{m, \ell}(  \Vec{x}^{\ell, k-1}  ) - \Vec{b}^\ell \right), \qquad  \text{ for } k = 1, 2, \dots, n_s^\ell.
\end{equation}
Let 
\[
\Vec{r}^{m, \ell, k} 
= 
\begin{bmatrix}
\Vec{r}_\sigma^{m, \ell, k}   \\
\Vec{r}_p^{m, \ell, k}  \\
\Vec{r}_s^{m, \ell, k}
\end{bmatrix}
:=
\Vec{r}^{m, \ell}(  \Vec{x}^{\ell, k-1} ) - \Vec{b}^\ell.
\]
The Jacobian system to be solved in each Newton step has the form:
\begin{equation}
\begin{split}
&\partial \Vec{r}^{m, \ell}(\Vec{x}^{\ell, k-1})\Delta \Vec{x}^{\ell, k}
=
\begin{bmatrix}
\frac{\partial \Vec{r}_\sigma^{m, \ell}}{\partial \Vec{\sigma}} &  \chak{-D^T}  & \frac{\partial \Vec{r}_\sigma^{m, \ell}}{\partial \Vec{s}} \\
D & 0 & 0 \\
\frac{\partial \Vec{r}_s^{m, \ell}}{\partial \Vec{\sigma}} & 0 & \frac{\partial \Vec{r}_s^{m, \ell}}{\partial \Vec{s}}
\end{bmatrix}
\begin{bmatrix}
\Delta \Vec{\sigma}^{\ell, k}  \\
\Delta \Vec{p}^{\ell, k}  \\
\Delta \Vec{s}^{\ell, k} 
\end{bmatrix}
= \begin{bmatrix}
\Vec{r}_\sigma^{m, \ell, k}   \\
\Vec{r}_p^{m, \ell, k}  \\
\Vec{r}_s^{m, \ell, k}
\end{bmatrix}
\end{split}
\label{eq:jacobian}
\end{equation}
The construction of each of the components of the Jacobian system is given in \ref{app:jacobian}.

\subsection{Iterative linear solver for the Jacobian system}

Because of the difference of $\frac{\partial \Vec{r}_\sigma^{m, \ell}}{\partial \Vec{\sigma}}$ on the fine and coarse levels, we will use two different linear solvers for the Jacobian system \eqref{eq:jacobian}, depending on the level $\ell$.
In both cases, \eqref{eq:jacobian} is transformed into a system involving only pressure (or face pressure) and saturation.
Then, a CPR-type preconditioner \cite{wallis1983incomplete,wallis1985constrained,lacroix2001decoupling,scheichl2003decoupling,cao2005} is applied to solve the transformed system.

\subsubsection{Transformation on the fine level}
When $\ell = 0$, $\frac{\partial \Vec{r}_\sigma^{m, \ell}}{\partial \Vec{\sigma}}$ is diagonal, so it is inexpensive to first eliminate $\Delta \Vec{\sigma}^{\ell, k}$ in \eqref{eq:jacobian} and obtain
\begin{equation}
\begin{split}
\blkMat{A}^\ell
\begin{bmatrix}
\Delta \Vec{p}^{\ell, k} \\
\Delta \Vec{s}^{\ell, k} 
\end{bmatrix}
:=
\begin{bmatrix}
A^\ell_{11} & A^\ell_{12} \\
A^\ell_{21} & A^\ell_{22}
\end{bmatrix}
\begin{bmatrix}
\Delta \Vec{p}^{\ell, k} \\
\Delta \Vec{s}^{\ell, k}
\end{bmatrix}
=
\begin{bmatrix}
\Vec{\xi}_1^{\ell, k} \\
\Vec{\xi}_2^{\ell, k}
\end{bmatrix}
\end{split}
\label{eq:primal_reduced}
\end{equation}
where
\begin{equation*}
\begin{split}
\begin{bmatrix}
A^\ell_{11} & A^\ell_{12} \\
A^\ell_{21} & A^\ell_{22}
\end{bmatrix}
:=
\begin{bmatrix}
0 & 0  \\
0 & \frac{\partial \Vec{r}_s^{m, \ell}}{\partial \Vec{s}}
\end{bmatrix}
-
\begin{bmatrix}
D \\
\frac{\partial \Vec{r}_s^{m, \ell}}{\partial \Vec{\sigma}} 
\end{bmatrix}
\left( \frac{\partial \Vec{r}_\sigma^{m, \ell}}{\partial \Vec{\sigma}} \right)^{-1} 
\begin{bmatrix}
\chak{-D^T} & \frac{\partial \Vec{r}_\sigma^{m, \ell}}{\partial \Vec{s}} 
\end{bmatrix},
\end{split}
\label{eq:hybrid_reduced_system}
\end{equation*}
and
\begin{equation*}
\begin{split}
\begin{bmatrix}
\Vec{\xi}_1^{\ell, k} \\
\Vec{\xi}_2^{\ell, k}
\end{bmatrix}
:=
\begin{bmatrix}
\Vec{r}_p^{m, \ell, k}  \\
\Vec{r}_s^{m, \ell, k} 
\end{bmatrix}
-
\begin{bmatrix}
D \\
\frac{\partial \Vec{r}_s^{m, \ell}}{\partial \Vec{\sigma}} 
\end{bmatrix}
\left( \frac{\partial \Vec{r}_\sigma^{m, \ell}}{\partial \Vec{\sigma}} \right)^{-1} 
\Vec{r}_\sigma^{m, \ell, k} \\
\end{split}
\label{eq:hybrid_reduced_rhs}
\end{equation*}

\subsubsection{Transformation on the coarse level}

When $\ell \ge 1$, $\frac{\partial \Vec{r}_\sigma^{m, \ell}}{\partial \Vec{\sigma}}$ is not a diagonal matrix.
Inverting $\frac{\partial \Vec{r}_\sigma^{m, \ell}}{\partial \Vec{\sigma}}$ would be too expensive.
Instead, we can consider the hybrid version of \eqref{eq:jacobian}, where the flux is replaced by one-sided fluxes.
Each one-sided flux is associated with one cell only, and the weak continuity of the flux is enforced through Lagrange multiplier (face pressure).
In the hybrid formulation, the system to be solved, conceptually, is
\begin{equation}
\begin{split}
\begin{bmatrix}
\widehat{\frac{\partial \Vec{r}_\sigma^m}{\partial \Vec{\sigma}}} &  \chak{-\widehat{D}^T}  & C^T & \widehat{\frac{\partial \Vec{r}_\sigma^m}{\partial \Vec{s}}} \\
\widehat{D} & 0 & 0 & 0 \\
C & 0 & 0 & 0 \\
\widehat{\frac{\partial \Vec{r}_s^m}{\partial \Vec{\sigma}}} & 0 & 0 & \frac{\partial \Vec{r}_s^m}{\partial \Vec{s}}
\end{bmatrix}
\begin{bmatrix}
\Delta \Vec{\sigma}^{\ell, k}  \\
\Delta \Vec{p}^{\ell, k}  \\
\Delta \Vec{\lambda}^{\ell, k} \\
\Delta \Vec{s}^{\ell, k}
\end{bmatrix}
& =
\begin{bmatrix}
\widehat{\Vec{r}_\sigma^{m, \ell, k}}  \\
\Vec{r}_p^{m, \ell, k} \\
\Vec{0} \\
\Vec{r}_s^{m, \ell, k}
\end{bmatrix}
\end{split}
\label{eq:hybrid}
\end{equation}
Since $\begin{bmatrix}
\widehat{\frac{\partial \Vec{r}_\sigma^m}{\partial \Vec{\sigma}}}  &  \chak{-\widehat{D}^T}  \\
\widehat{D} & 0 \\
\end{bmatrix}$ is block-diagonal and each block is invertible, the global system to be solved in our implementation is the reduced system
\begin{equation}
\begin{split}
\blkMat{A}^\ell
\begin{bmatrix}
\Delta \Vec{\lambda}^{m, k} \\
\Delta \Vec{s}^{m, k}
\end{bmatrix}
:=
\begin{bmatrix}
A^\ell_{11} & A^\ell_{12} \\
A^\ell_{21} & A^\ell_{22}
\end{bmatrix}
\begin{bmatrix}
\Delta \Vec{\lambda}^{\ell, k}\\
\Delta \Vec{s}^{\ell, k} 
\end{bmatrix}
=
\begin{bmatrix}
\Vec{\xi}_1^{\ell, k} \\
\Vec{\xi}_2^{\ell, k}
\end{bmatrix}
\end{split}
\label{eq:hybrid_reduced}
\end{equation}
where
\begin{equation*}
\begin{split}
\begin{bmatrix}
A^\ell_{11} & A^\ell_{12} \\
A^\ell_{21} & A^\ell_{22}
\end{bmatrix}
:=
\begin{bmatrix}
0 & 0 \\
0 & \frac{\partial \Vec{r}_s^m}{\partial \Vec{s}}
\end{bmatrix}
-
\begin{bmatrix}
C & 0 \\
\widehat{\frac{\partial \Vec{r}_s^m}{\partial \Vec{\sigma}}} & 0
\end{bmatrix}
\begin{bmatrix}
\widehat{\frac{\partial \Vec{r}_\sigma^m}{\partial \Vec{\sigma}}}  &  \chak{-\widehat{D}^T}  \\
\widehat{D} & 0 \\
\end{bmatrix}^{-1}
\begin{bmatrix}
C^T & \widehat{\frac{\partial \Vec{r}_\sigma^m}{\partial \Vec{s}}} \\
 0 & 0 \\
\end{bmatrix},
\end{split}
\label{eq:hybrid_reduced_system}
\end{equation*}
and
\begin{equation*}
\begin{split}
\begin{bmatrix}
\Vec{\xi}_1^{\ell, k} \\
\Vec{\xi}_2^{\ell, k}
\end{bmatrix}
:=
\begin{bmatrix}
\Vec{r}_s^{m, \ell, k}  \\
\Vec{0}
\end{bmatrix}
-
\begin{bmatrix}
C & 0 \\
\widehat{\frac{\partial \Vec{r}_s^m}{\partial \Vec{\sigma}}} & 0 
\end{bmatrix}
\begin{bmatrix}
\widehat{\frac{\partial \Vec{r}_\sigma^m}{\partial \Vec{\sigma}}}  &  \chak{-\widehat{D}^T}  \\
\widehat{D} & 0 \\
\end{bmatrix}^{-1}
\begin{bmatrix}
\widehat{\Vec{r}_\sigma^{m, \ell, k}} \\
\Vec{r}_p^{m, \ell, k} \\
\end{bmatrix}.
\end{split}
\label{eq:hybrid_reduced_rhs}
\end{equation*}

\subsubsection{Linear solver for the transformed system}

Systems \eqref{eq:primal_reduced} and \eqref{eq:hybrid_reduced} are solved by preconditioned GMRES.
Since \eqref{eq:primal_reduced} (respectively \eqref{eq:hybrid_reduced} is a system involving pressure (respectively face pressure) and saturation, a CPR-type preconditioner \cite{wallis1983incomplete,wallis1985constrained,lacroix2001decoupling,scheichl2003decoupling,cao2005} is employed. 
More precisely, we consider a two-stage preconditioner
\begin{equation}
\blkMat{B}^\ell := \blkMat{B}^\ell_1 + \blkMat{B}^\ell_2 (\blkMat{I} - \blkMat{A}^\ell \blkMat{B}^\ell_1)
\end{equation}
which gives rise to the product iteration matrix
\begin{equation*}
    (\blkMat{I}-\blkMat{A}^\ell\blkMat{B}^\ell)=
    (\blkMat{I}-\blkMat{A}^\ell\blkMat{B}^\ell_2)(\blkMat{I}-\blkMat{A}^\ell\blkMat{B}^\ell_1).
\end{equation*}
The first stage is a block lower-triangular preconditioner
\begin{equation*}
\blkMat{B}^\ell_1 := \begin{bmatrix}
B^\ell_{11} &   0 \\
B^\ell_{22}A^\ell_{12}B^\ell_{11} & B^\ell_{22}
\end{bmatrix},
\end{equation*}
where $B^\ell_{11}$ is an AMG preconditioner \cite{hypre,stuben2007algebraic} for $A^\ell_{11}$, and $B^\ell_{22}$ is the $\ell_1$-Jacobi smoother \cite{baker11} for $A^\ell_{22}$.
The second stage is an ILU(1) preconditioner for the monolithic system $\blkMat{A}^\ell$. That is, $\blkMat{B}^\ell_2 := ILU(1)(\blkMat{A}^\ell) $.

\section{Numerical examples}\label{sec:numerics}

In this section, we consider three challenging test cases to demonstrate the performance of the FAS-based nonlinear multigrid algorithm.
The test cases are selected to illustrate the applicability of the method to realistic reservoir simulation problems, as well as its robustness with respect to specific numerical challenges often encountered in subsurface flow and transport applications.
In Section~\ref{sec:spe10_bottom_layer}, we use layer 85 of the SPE10 test case \cite{spe10} to demonstrate the behavior of the scheme for a highly heterogeneous geological model.
We also consider two distinct sets of fluid parameters to show that the nonlinear algorithm is robust in different mobility regimes and can handle the propagation of sharp saturation fronts in the domain.
In Section~\ref{sec:egg_model}, we use the Egg model \cite{EggModel} to illustrate that FAS achieves excellent performance for a wide range of relative permeability parameters controlling the nonlinearity in the problem.
Section~\ref{sec:saigup_model} highlights the ability of the method to handle the high geometric complexity inherent in corner-point grids using the SAIGUP model \cite{manzocchi2008sensitivity}.
The parameters employed in the simulations and the problem sizes can be found in Tables~\ref{tab:parameters} and \ref{tab:problem_sizes}, respectively.

Even though the fine-scale meshes used in the three test cases are structured, coarse cell aggregates are generated using METIS \cite{karypis1998fast} without relying on any intrinsic structure.
Before calling METIS, we remove the cells that are connected to the wells from the cell-connectivity graph (see Remark~\ref{remark:wells}).
The coarsening factor $\beta$ reported is computed using the average aggregate size on each level.
To make the test cases challenging, we use an aggressive time stepping strategy in which the time step size is multiplied by a factor $\nu >1$ at every step:
\begin{equation}
\Delta t_{m} = \nu \Delta t_{m-1}, \quad m \geq 1. \label{time_stepping}
\end{equation}
The corresponding Courant-Friedrichs-Lewy (CFL) numbers are computed using the standard formula \cite{cao2002development}.
Nonlinear convergence is achieved when the normalized residual drops below $10^{-6}$.
The maximum number of nonlinear iterations is set to 10 on the coarsest level, and to 1 on all other levels.
To enhance nonlinear convergence, we force the saturations to remain in [0,1] after each fine-level FAS update (referred to as local saturation chopping, see \cite{younis2011modern}).
After the coarse-level FAS updates, we do not use local saturation chopping and instead we extend the mobility functions with constant values outside [0,1].
The FAS results are compared to those obtained with single-level Newton with local saturation chopping. 
%

The discrete problems are generated using our own implementation based on MFEM \cite{mfem} of the finite-volume scheme described in Section \ref{sec:FV_tpfa}.
The multilevel spectral coarsening is performed with smoothG \cite{smoothg}, and the visualization is generated with GLVis \cite{glvis}.

\begin{table}[htbp]
	\caption{Parameter values used for the numerical examples. For the SPE10 test, we consider a case with a favorable end-point mobility ratio, and a case with an unfavorable end-point mobility ratio (see the definition in Section~\ref{sec:spe10_bottom_layer}). For the three test cases, the times are reported in total pore volume injected (PVI), which is the ratio of the injected wetting-phase volume over the total pore volume of the reservoir.}
	\label{tab:parameters}
        \small
        \centering
        \begin{tabular}{lllllllll}
          \toprule
          Symbol & Parameter & Units &
          \begin{tabular}{@{}c@{}}SPE10 layer 85 \\ (favorable) \end{tabular} & 
          \begin{tabular}{@{}c@{}}SPE10 layer 85 \\ (unfavorable) \end{tabular} &         
          Egg & 
          SAIGUP\\ 
          \midrule
          $s^0$        &  Initial wetting-phase saturation    & [-]           & 0 & 0 & 0 & 0 \\
          $\mu_w$      &  Wetting-phase viscosity             & [Pa.s]        & $10^{-3}$ & $10^{-3}$ & $10^{-3}$ & $10^{-3}$ \\
          $\mu_{nw}$   &  Non-wetting phase viscosity         & [Pa.s]        & $2.0 \times 10^{-4}$ & $5.0 \times 10^{-3}$ & $5.0 \times 10^{-3}$ & $5.0 \times 10^{-3}$ \\          
          $\gamma$     &  Relative permeability exponent      & [-]           & 2 & 2 & 2, 3, or 4 & 2 \\
          $\lambda_{\alpha}$&  Phase mobility ($\alpha \in \{w,nw\}$)         & [Pa$^{-1}$.s$^{-1}$] &  $s^{\gamma}_{\alpha}/\mu_{\alpha}$ & $s^{\gamma}_{\alpha}/\mu_{\alpha}$ &  $s^{\gamma}_{\alpha}/\mu_{\alpha}$ & $s^{\gamma}_{\alpha}/\mu_{\alpha}$ \\
          \midrule 
          $n^I$        &  Number of injectors                 & [-]           & 1 & 1 & 8 & 5 \\
          $n^P$        &  Number of producers                 & [-]           & 4 & 4 & 4 & 5 \\
          $q_w^I$      &  Wetting-phase injection rate        & [m$^3$.s$^{-1}$] & $1.8 \times 10^{-4}$ & $6.1 \times 10^{-5}$ & 10$^{-3}$ & $10^{-1}$ \\
          $p_{bh}$     &  Bottomhole pressure                 & [Pa]          & $10^6$ & $10^6$ & $10^6$ & $10^6$ \\
          \midrule
          $\Delta t_0$ &  Initial time step                   & [PVI]         & $5.3 \times 10^{-4}$ & $1.8 \times 10^{-4}$ & $9.1 \times 10^{-5}$ & $4.2 \times 10^{-5}$ \\
          $\nu$          &  Time step increase factor           & [-]           & 2, 4, or 8 & 2, 4, or 8 & 2 & 2 \\
          $T_{f}$      &  Final time                          & [PVI]         & $2.7 \times 10^{-1}$ & $9.0 \times 10^{-2}$ & $4.6 \times 10^{-2}$ & $2.1 \times 10^{-2}$ \\
          \bottomrule
        \end{tabular}
\end{table}


%

\begin{table}[htbp]
	\caption{Problem sizes in the numerical examples.}
	\label{tab:problem_sizes}
        \small
        \centering
        \begin{tabular}{lrrrrrr}
          \toprule
          &
          \begin{tabular}{@{}c@{}}SPE10 \\ layer 85 \end{tabular} &
          \begin{tabular}{@{}c@{}}Egg \\ (refined) 
          \end{tabular} &          
          \begin{tabular}{@{}c@{}}SAIGUP \\ (refined) \end{tabular}\\
          \midrule
          $|\mathcal{T}|$               & 12,321 & 148,424 & 629,760  \\
          $|\mathcal{E}|$               & 25,355 & 431,092 & 1,912,471 \\
          Number of unknowns            & 37,676 & 579,516 & 2,542,231 \\
          \bottomrule          
        \end{tabular}
\end{table}


\subsection{SPE10 layer 85}\label{sec:spe10_bottom_layer}

\begin{figure}[ht]
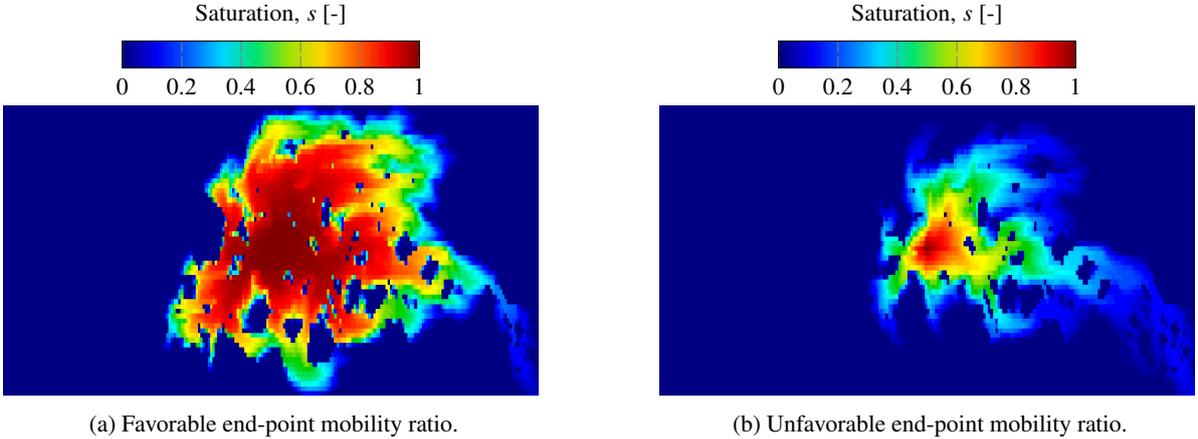

  \small
  \centering
  \begin{subfigure}[b]{0.475\textwidth}
    \centering
    \if \generateTikzFigures 1
      \include{./pics/paper_pics/sat-spe-bottom-sharp}    
    \else
      \includegraphics[scale=1]{main-figure\theFigureCounter.pdf}
      \stepcounter{FigureCounter} 
    \fi
    \caption{{Favorable end-point mobility ratio.}}
  \end{subfigure}
  \hfill
  \begin{subfigure}[b]{0.475\textwidth}
    \centering
    \if \generateTikzFigures 1    
      \include{./pics/paper_pics/sat-spe-bottom-smeared}  
    \else
      \includegraphics[scale=1]{main-figure\theFigureCounter.pdf}
      \stepcounter{FigureCounter}  
    \fi
    \caption{Unfavorable end-point mobility ratio.}
  \end{subfigure}
\caption{Final wetting-phase saturation maps in layer 85 of the SPE10 test case.}
\label{fig:sat-spe-bottom}
\end{figure}

\begin{figure}[h!]
  \small
  \centering
  \begin{subfigure}[b]{0.475\textwidth}
    \centering
    \if \generateTikzFigures 1    
      \include{./pics/paper_pics/spe_sharp_CFL_PVI}    
    \else
      \includegraphics[scale=1]{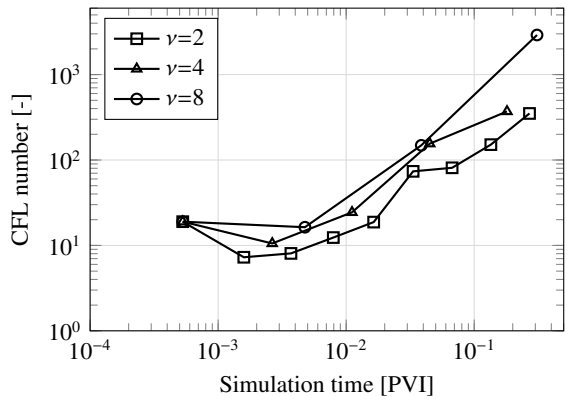}
      \stepcounter{FigureCounter} 
    \fi
    \caption{Favorable end-point mobility ratio.}
  \end{subfigure}
  \hfill
  \begin{subfigure}[b]{0.475\textwidth}
    \centering
    \if \generateTikzFigures 1
      \include{./pics/paper_pics/spe_smeared_CFL_PVI}  
    \else
      \includegraphics[scale=1]{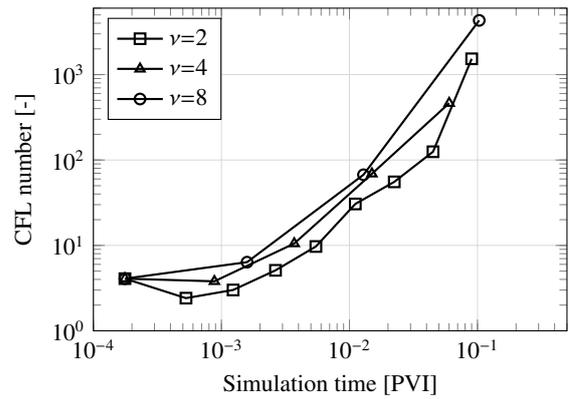}
      \stepcounter{FigureCounter}  
    \fi
    \caption{Unfavorable end-point mobility ratio.}
  \end{subfigure}
  \caption{CFL numbers [-] as a function of simulation time [PVI] for the three time stepping strategies in the SPE10 test case (layer 85). As $\nu$ increases, the time stepping strategy becomes more aggressive (see Eq.~\eqref{time_stepping}).}
  \label{fig:spe10_bot_CFL_compare_semilogx}
\end{figure}

We first consider layer 85 of the SPE10 model \cite{spe10} to assess the performance of FAS on highly heterogeneous permeability and porosity fields.
Following standard practice, the cells with a pore volume smaller than a threshold of 0.5 m$^3$ are treated as inactive and are removed from the mesh.
The wells are placed according to the specifications of the original test case, with an injector at the center and a producer in each corner of the domain.
The well Peaceman indices are computed with MRST \cite{krogstad2015mrst}. 
The propagation of sharp saturation fronts is a well-documented challenge for multilevel solution algorithms applied to multiphase flows \cite{aarnes2004use,kippe2008comparison}.
To illustrate the robustness of FAS for various flow regimes, we consider two sets of fluid parameters (see Table~\ref{tab:parameters}), named using the classical petroleum engineering terminology:
\begin{itemize}
\item Favorable end-point mobility ratio: the non-wetting phase viscosity is set to $\mu_{nw} = 2 \times 10^{-4}$ Pa.s. %
  The resulting end-point mobility ratio of  $\lambda(1)/\lambda(0) = 5$ leads to a flow regime characterized by the propagation of sharp, piston-like saturation fronts. %
\item Unfavorable end-point mobility ratio: the non-wetting phase viscosity is set to $\mu_{nw} = 5 \times 10^{-3}$ Pa.s. %
  This choice yields an end-point mobility ratio of $\lambda(1)/\lambda(0) = 0.2$, producing a flow regime characterized by the propagation of smeared fronts with small-scale saturation fingers. %
\end{itemize}
To obtain similar CFL numbers with the two sets of parameters, we use two different injection rates (see Table~\ref{tab:parameters}). 
The final saturation maps for the two scenarios are presented in Fig.~\ref{fig:sat-spe-bottom}.
To show that FAS can handle very large time step sizes, we multiply the time step size by $\nu \in \{2, 4, 8\}$ at the end of each step, starting with a time step 
corresponding to a CFL number of approximately 4 and 20 for the favorable and unfavorable end-point mobility ratio cases, respectively.
The simulation involves nine, five, and four time steps for $\nu = 2$, 4, and 8, respectively.
These aggressive time stepping strategies, combined with quadratic phase mobilities, produce an increasingly difficult test case for the nonlinear solvers.
Although the resulting CFL numbers may seem very large toward the end of the simulation (see Fig.~\ref{fig:spe10_bot_CFL_compare_semilogx}), our goal is to check that FAS converges quickly regardless of the time step size.  With such a method, the time stepping strategy can then be decided only based on accuracy considerations--i.e. solely focused on limiting temporal truncation errors, rather than improving nonlinear convergence.


\begin{figure}[h!]
  \small
  \centering
  \begin{subfigure}[b]{0.32\textwidth}
    \centering
    \if \generateTikzFigures 1
      \include{./pics/paper_pics/spe_sharp_k2_NLIter}    
     \else
      \includegraphics[scale=1]{main-figure\theFigureCounter.pdf}
      \stepcounter{FigureCounter}
    \fi
    \caption{$\nu=2$.}
  \end{subfigure}
  \hfill
  \begin{subfigure}[b]{0.32\textwidth}
    \centering
    \if \generateTikzFigures 1
      \include{./pics/paper_pics/spe_sharp_k4_NLIter}   
    \else
      \includegraphics[scale=1]{main-figure\theFigureCounter.pdf}
      \stepcounter{FigureCounter}   
    \fi
    \caption{$\nu=4$.}
  \end{subfigure}
  \hfill
  \begin{subfigure}[b]{0.32\textwidth}
    \centering
    \if \generateTikzFigures 1
      \include{./pics/paper_pics/spe_sharp_k8_NLIter} 
    \else
      \includegraphics[scale=1]{main-figure\theFigureCounter.pdf}
      \stepcounter{FigureCounter}   
    \fi
    \caption{$\nu=8$.}
  \end{subfigure}
  \caption{Number of nonlinear iterations per time step as a function of CFL number [-] for the favorable end-point mobility ratio in the SPE10 test case (layer 85).}
  \label{fig:spe10_bot_CFL_iterations_sharp}
\end{figure}

\begin{figure}[h!]
  \small
  \centering
  \begin{subfigure}[b]{0.32\textwidth}
    \centering
    \if \generateTikzFigures 1
      \include{./pics/paper_pics/spe_smear_k2_NLIter}    
    \else
      \includegraphics[scale=1]{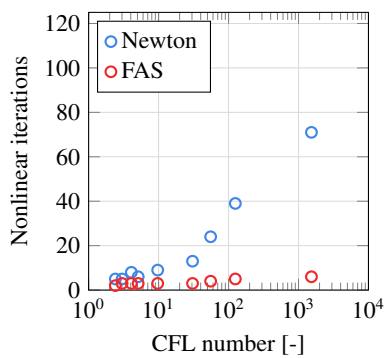}
      \stepcounter{FigureCounter}
    \fi
    \caption{$\nu=2$.}
  \end{subfigure}
  \hfill
  \begin{subfigure}[b]{0.32\textwidth}
    \centering
    \if \generateTikzFigures 1
      \include{./pics/paper_pics/spe_smear_k4_NLIter}  
    \else
      \includegraphics[scale=1]{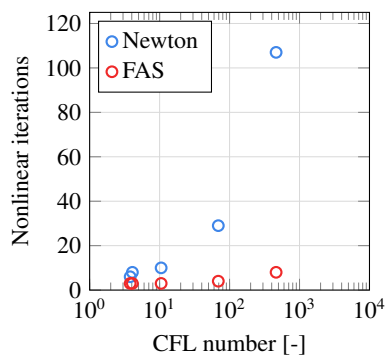}
      \stepcounter{FigureCounter}   
    \fi
    \caption{$\nu=4$.}
  \end{subfigure}
  \hfill
  \begin{subfigure}[b]{0.32\textwidth}
    \centering
    \if \generateTikzFigures 1
      \include{./pics/paper_pics/spe_smear_k8_NLIter} 
    \else
      \includegraphics[scale=1]{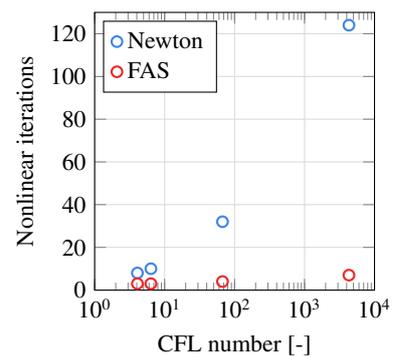}
      \stepcounter{FigureCounter}   
    \fi
    \caption{$\nu=8$.}
  \end{subfigure}
  \caption{Number of nonlinear iterations per time step as a function of CFL number [-] for the unfavorable end-point mobility ratio in the SPE10 test case (layer 85).}
  \label{fig:spe10_bot_CFL_iterations_smear}
\end{figure}

In Figs.~\ref{fig:spe10_bot_CFL_iterations_sharp} and \ref{fig:spe10_bot_CFL_iterations_smear}, we compare the number of nonlinear iterations per step required by FAS and single-level Newton.
FAS is configured with two levels and a coarsening factor $\beta = 16$.
We observe that as the time step becomes larger, with CFL numbers larger than 20, the number of nonlinear iterations per step required by single-level Newton increases drastically for both the favorable and unfavorable end-point mobility ratios.
Using the most aggressive time stepping strategy ($\nu = 8$), single-level Newton needs 139 iterations and 124 iterations to reach convergence for the last time step in the favorable and unfavorable cases, respectively. 
However, with FAS, the nonlinear behavior is more stable throughout the simulation.  It only exhibits a moderate increase in the number of nonlinear iterations for the challenging favorable end-point mobility ratio case performed with $\nu = 8$, reaching 13 iterations for the last time step.
In the other configurations, the number of FAS iterations remains less than 8 for all time steps.

\begin{figure}[h!]
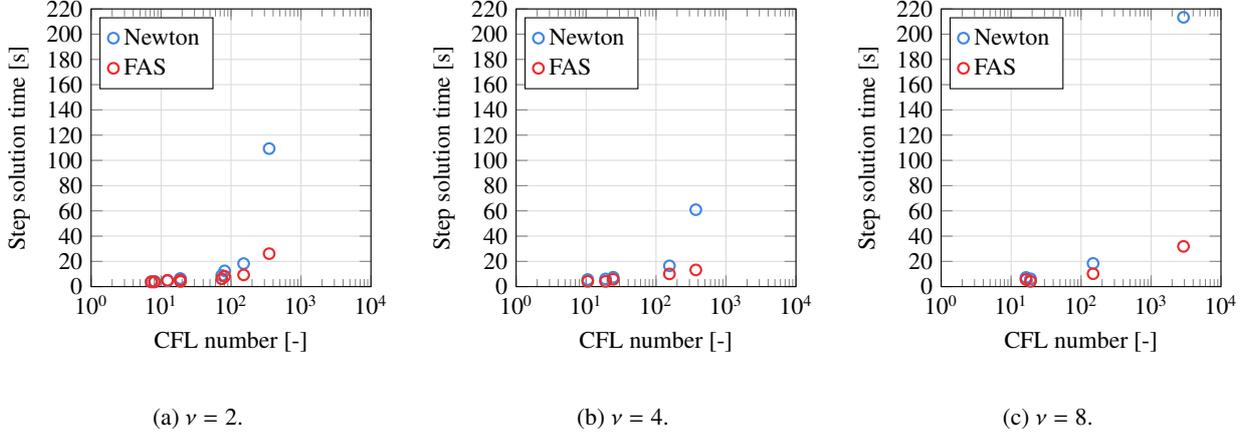

  \small
  \centering
  \begin{subfigure}[b]{0.32\textwidth}
    \centering
    \if \generateTikzFigures 1
      \include{./pics/paper_pics/spe_sharp_k2_tSolv}    
    \else
      \includegraphics[scale=1]{main-figure\theFigureCounter.pdf}
      \stepcounter{FigureCounter} 
    \fi
    \caption{$\nu=2$.}
  \end{subfigure}
  \hfill
  \begin{subfigure}[b]{0.32\textwidth}
    \centering
    \if \generateTikzFigures 1
      \include{./pics/paper_pics/spe_sharp_k4_tSolv}    
    \else
      \includegraphics[scale=1]{main-figure\theFigureCounter.pdf}
      \stepcounter{FigureCounter}   
    \fi
    \caption{$\nu=4$.}
  \end{subfigure}
  \hfill
  \begin{subfigure}[b]{0.32\textwidth}
    \centering
    \if \generateTikzFigures 1
      \include{./pics/paper_pics/spe_sharp_k8_tSolv}    
    \else            
      \includegraphics[scale=1]{main-figure\theFigureCounter.pdf}
      \stepcounter{FigureCounter}   
    \fi
    \caption{$\nu=8$.}
  \end{subfigure}
  \caption{Step solution time [s] as a function of CFL number [-] for the favorable end-point mobility ratio in the SPE10 test case (layer 85).}
  \label{fig:spe10_bot_CFL_time_sharp}
\end{figure}

In Figs.~\ref{fig:spe10_bot_CFL_time_sharp} and \ref{fig:spe10_bot_CFL_time_smear}, we consider the solution time per step for the two solution algorithms.
We note that for CFL numbers smaller than 20, the solution time per step of FAS and single-level Newton is similar.
However, for CFL numbers larger than 20, the increase in nonlinear iterations observed in Figs.~\ref{fig:spe10_bot_CFL_iterations_sharp}-\ref{fig:spe10_bot_CFL_iterations_smear} for single-level Newton produces a sharp rise in solution time per step.
This is not the case for FAS, whose robust nonlinear behavior for large CFL numbers limits the increase in solution time per step.
As a result, FAS yields a significant reduction in the solution time per step compared to single-level Newton for these large time steps.
The last time step of the simulation performed with $\nu=8$ provides a good illustration of the superior behavior of the multigrid algorithm, since FAS is 6.7 times and 12.0 times faster than single-level Newton for the favorable and unfavorable cases, respectively.

\begin{figure}[h!]
  \small
  \centering
  \begin{subfigure}[b]{0.475\textwidth}
    \centering
    \if \generateTikzFigures 1
      \include{./pics/paper_pics/spe_sharp_CFL_tSolvCumul}    
    \else
      \includegraphics[scale=1]{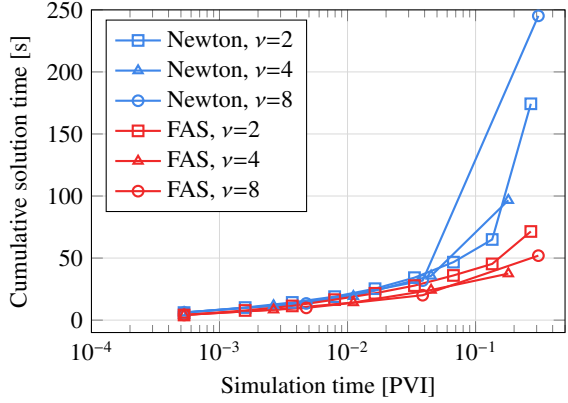}
      \stepcounter{FigureCounter} 
    \fi
    \caption{Favorable end-point mobility ratio.}
  \end{subfigure}
  \hfill
  \begin{subfigure}[b]{0.475\textwidth}
    \centering
    \if \generateTikzFigures 1
      \include{./pics/paper_pics/spe_smear_CFL_tSolvCumul}  
    \else
      \includegraphics[scale=1]{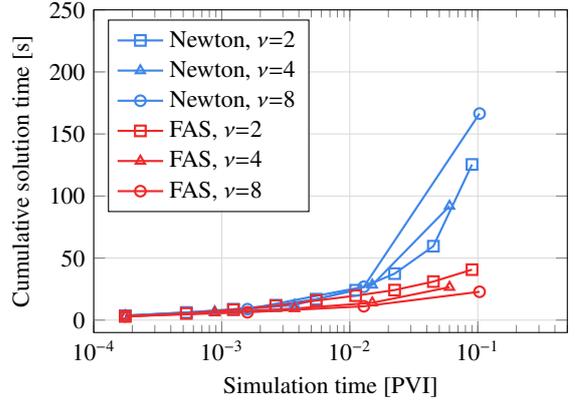}
      \stepcounter{FigureCounter}   
    \fi
    \caption{Unfavorable end-point mobility ratio.}
  \end{subfigure}
  \caption{Cumulative solution time [s] as a function of simulation time [PVI] for FAS and single-level Newton in the SPE10 test case (layer 85).}
  \label{fig:spe10_bot_k_compare_semilogx}
\end{figure}

Figure~\ref{fig:spe10_bot_k_compare_semilogx} summarizes the performance comparison between the two nonlinear solvers by showing the cumulative solution time as a function of simulation time.
Thanks to its robustness for large time steps---corresponding to large CFL numbers---FAS yields a significant reduction in the total solution time compared to single-level Newton for all values of $\nu$. 
This conclusion holds for the two flow regimes considered in our study.

\begin{figure}[h!]
  \small
  \centering
  \begin{subfigure}[b]{0.32\textwidth}
    \centering
    \if \generateTikzFigures 1
      \include{./pics/paper_pics/spe_smear_k2_tSolv}    
    \else
      \includegraphics[scale=1]{main-figure\theFigureCounter.pdf}
      \stepcounter{FigureCounter}
    \fi
    \caption{$\nu=2$.}
  \end{subfigure}
  \hfill
  \begin{subfigure}[b]{0.32\textwidth}
    \centering
    \if \generateTikzFigures 1
      \include{./pics/paper_pics/spe_smear_k4_tSolv}    
    \else
      \includegraphics[scale=1]{main-figure\theFigureCounter.pdf}
      \stepcounter{FigureCounter}   
    \fi
    \caption{$\nu=4$.}
  \end{subfigure}
  \hfill
  \begin{subfigure}[b]{0.32\textwidth}
    \centering
    \if \generateTikzFigures 1
      \include{./pics/paper_pics/spe_smear_k8_tSolv}    
    \else
      \includegraphics[scale=1]{main-figure\theFigureCounter.pdf}
      \stepcounter{FigureCounter}   
    \fi
    \caption{$\nu=8$.}
  \end{subfigure}
  \caption{Step solution time [s] as a function of CFL number [-] for the unfavorable end-point mobility ratio in the SPE10 test case (layer 85).}
  \label{fig:spe10_bot_CFL_time_smear}
\end{figure}

\subsection{The Egg model}\label{sec:egg_model}

\begin{figure}
 \small
 \centering
 \begin{subfigure}[b]{0.3\textwidth}
   \centering    
   \if \generateTikzFigures 1
     \include{./pics/paper_pics/egg_final_saturation_gamma4}
   \else
      \includegraphics[scale=1]{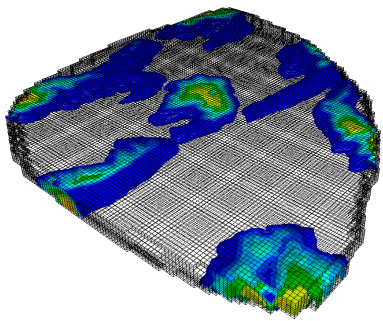}
      \stepcounter{FigureCounter} 
   \fi
   \caption{Refined Egg model.}
   \label{fig:sat-egg-smeared}
 \end{subfigure}
 \hfill
 \begin{subfigure}[b]{0.6\textwidth}
   \centering
   \if \generateTikzFigures 1
     \include{./pics/paper_pics/new_SAIGUP_fine_final_saturation_gamma2}
   \else
      \includegraphics[scale=1]{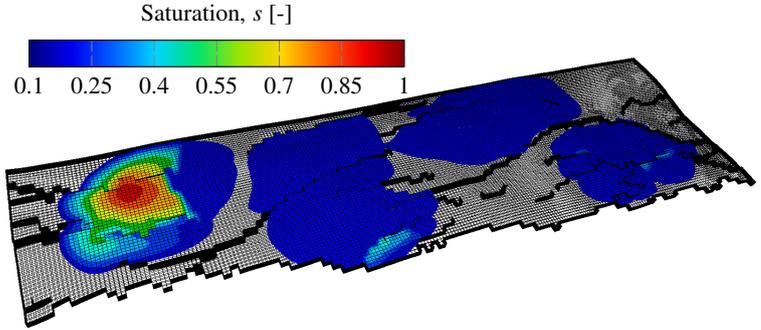}
      \stepcounter{FigureCounter} 
   \fi
   \caption{Refined SAIGUP model.}
   \label{fig:sat-saigup-smeared}
 \end{subfigure}
\caption{Final wetting-phase saturation field for the refined Egg (a) and SAIGUP (b) using quartic ($\gamma = 4$) and quadratic ($\gamma = 2$) relative permeabilities, respectively. 
}
\end{figure}

The topological and geological properties used in this section are derived from the Egg model \cite{EggModel}.
The simulations are performed on a refined  mesh consisting of 148,424 active cells generated with a 2 $\times$ 2 $\times$ 2 regular refinement of the original Egg model mesh consisting of 18,553 active cells (see Table~\ref{tab:problem_sizes}).
To increase the heterogeneity of the model, we rescale the permeability field and impose a ratio of $2 \times 10^5$ between the largest and smallest permeability values in each direction.
The porosity field is homogeneous with $\phi = 0.2$.
%
%
To evaluate the robustness of FAS when the strength of the nonlinearity increases, we consider three relative permeability exponents $\gamma \in \{ 2, 3, 4 \}$ in the analytical expression defining the phase mobility functions (see Table~\ref{tab:parameters}).
We select an unfavorable mobility ratio of 0.2, which produces a flow regime characterized by smeared saturation fronts. 
The wells are placed using the specifications of \cite{EggModel}, but we only keep one perforation per well, chosen as the perforation with the largest Peaceman index computed by MRST for each well.
We use the time stepping method of Eq.~\eqref{time_stepping}, with $\nu=2$. 
The final wetting-phase saturation map is presented in Fig.~\ref{fig:sat-egg-smeared}.

\begin{figure}[h!]
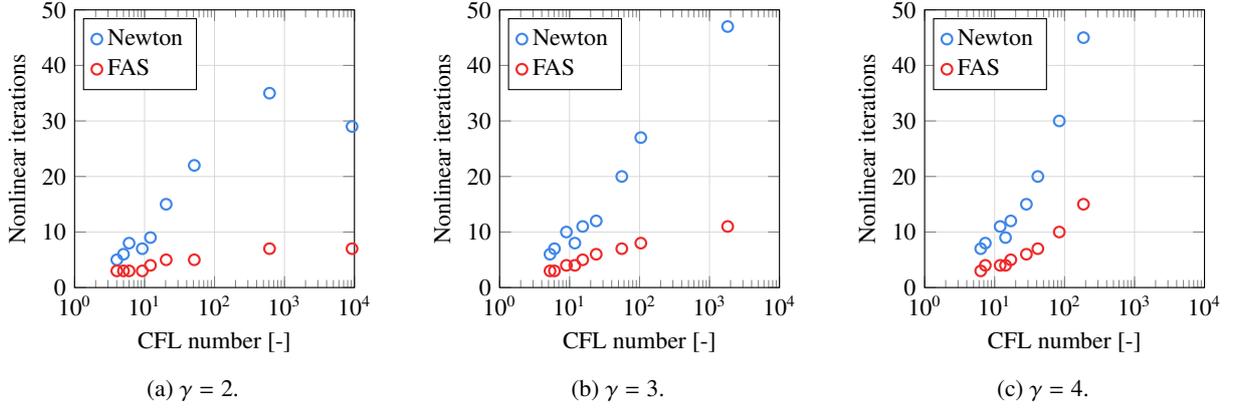

  \small
  \centering
  \begin{subfigure}[b]{0.32\textwidth}
    \centering
    \if \generateTikzFigures 1
      \include{./pics/paper_pics/egg_gamma2_NLIter}    
    \else
      \includegraphics[scale=1]{main-figure\theFigureCounter.pdf}
      \stepcounter{FigureCounter} 
    \fi
    \caption{$\gamma=2$.}
  \end{subfigure}
  \hfill
  \begin{subfigure}[b]{0.32\textwidth}
    \centering
    \if \generateTikzFigures 1
      \include{./pics/paper_pics/egg_gamma3_NLIter}  
    \else
      \includegraphics[scale=1]{main-figure\theFigureCounter.pdf}
      \stepcounter{FigureCounter}  
    \fi
    \caption{$\gamma=3$.}
  \end{subfigure}
  \hfill
  \begin{subfigure}[b]{0.32\textwidth}
    \centering
    \if \generateTikzFigures 1
      \include{./pics/paper_pics/egg_gamma4_NLIter} 
    \else
      \includegraphics[scale=1]{main-figure\theFigureCounter.pdf}
      \stepcounter{FigureCounter}  
    \fi
    \caption{$\gamma=4$.}
  \end{subfigure}
  \caption{Number of nonlinear iterations per time step as a function of CFL number [-] in the refined Egg test case.}
  \label{fig:egg_CFL_iterations}
\end{figure}

\begin{figure}[h!]
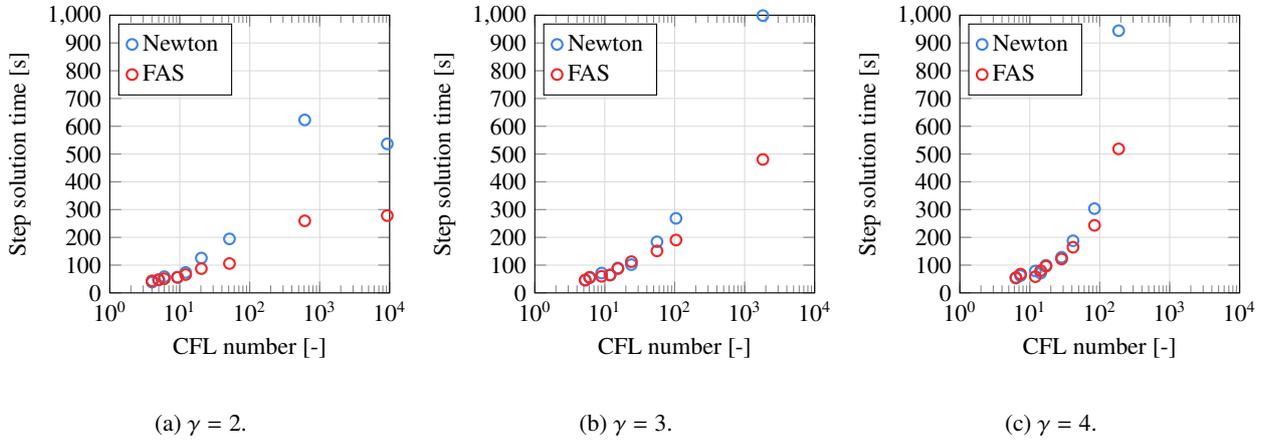

  \small
  \centering
  \begin{subfigure}[b]{0.32\textwidth}
    \centering
    \if \generateTikzFigures 1
      \include{./pics/paper_pics/egg_gamma2_tSolv}    
    \else
      \includegraphics[scale=1]{main-figure\theFigureCounter.pdf}
      \stepcounter{FigureCounter}
    \fi
    \caption{$\gamma=2$.}
  \end{subfigure}
  \hfill
  \begin{subfigure}[b]{0.32\textwidth}
    \centering
    \if \generateTikzFigures 1
      \include{./pics/paper_pics/egg_gamma3_tSolv}    
    \else
      \includegraphics[scale=1]{main-figure\theFigureCounter.pdf}
      \stepcounter{FigureCounter}   
    \fi
    \caption{$\gamma=3$.}
  \end{subfigure}
  \hfill
  \begin{subfigure}[b]{0.32\textwidth}
    \centering
    \if \generateTikzFigures 1
      \include{./pics/paper_pics/egg_gamma4_tSolv}    
    \else
      \includegraphics[scale=1]{main-figure\theFigureCounter.pdf}
      \stepcounter{FigureCounter}   
    \fi
    \caption{$\gamma=4$.}
  \end{subfigure}
  \caption{Step solution time [s] as a function of CFL number [-] in the refined Egg test case.}
  \label{fig:egg_CFL_time}
%

\end{figure}

We use a three-level FAS with a coarsening factor $\beta = 32$ for the refined mesh.
The nonlinear behavior and solve time per step for FAS and single-level Newton on the refined mesh are documented in Figs.~\ref{fig:egg_CFL_iterations} and \ref{fig:egg_CFL_time}.
The results are in agreement with those of the previous section.
As the time step size increases, FAS only exhibits a limited increase in the number of iterations, while the nonlinear behavior obtained with single-level Newton deteriorates quickly. 
%
For large CFL numbers, FAS achieves a large reduction in both nonlinear iteration counts and solving time for all values of $\gamma$ tested.
Table~\ref{tab:refined-egg} summarizes our observations. 
%
Compared with single-level Newton, FAS reduces the total solving time by respectively 43\%, 33\%, and 27\% for $\gamma = 2, 3$, and 4.


\begin{table}[!ht]
	\caption{Solution time, $T_{\text{sol}}$ [s], and average number of nonlinear iterations per time step ($n_{\text{it}}$) for the refined Egg model. FAS is based on a three-level hierarchy and a coarsening factor of $\beta = 32$.}
	\label{tab:refined-egg}
        \small
        \centering
        \begin{tabular}{lcccccc}
          \toprule
          Solver & \multicolumn{2}{c}{$\gamma = 2$} & \multicolumn{2}{c}{$\gamma = 3$} & \multicolumn{2}{c}{$\gamma = 4$} \\
          \cmidrule(l){2-3} \cmidrule(l){4-5} \cmidrule(l){6-7}
          & $T_{\text{sol}}$ & $n_{\text{it}}$ & $T_{\text{sol}}$ & $n_{\text{it}}$ & $T_{\text{sol}}$ & $n_{\text{it}}$ \\
          \midrule
          Single-level Newton & 1,753 & 15.11 & 1,878 & 16.44 & 1,928 & 17.44 \\
          FAS                 & 995 &  4.44 & 1,249 &  5.67 & 1,406 &  6.44 \\
          \bottomrule          
        \end{tabular}
\end{table}



\subsection{The SAIGUP model}\label{sec:saigup_model}

This numerical example is based on the SAIGUP model \cite{manzocchi2008sensitivity}.
It aims at demonstrating that FAS can handle the geometric complexity of realistic corner-point meshes widely used in industrial reservoir simulation studies.
In this example, we use a regularly refined version of the original mesh (consisting of 78,710 cells).
The refined mesh has a total of 629,760 cells.
The heterogeneous permeability and porosity fields are those of the original model.
The well placement is the same as in the original test case with only one perforation per well, computed by MRST and selected as in Section~\ref{sec:egg_model}.
We use the set of fluid properties yielding an unfavorable mobility ratio equal to 0.2.
We refer the reader to Table~\ref{tab:parameters} for the complete list of parameters used in this test case.
The final wetting-phase saturation field is presented in Fig.~\ref{fig:sat-saigup-smeared}.

\begin{figure}[ht]
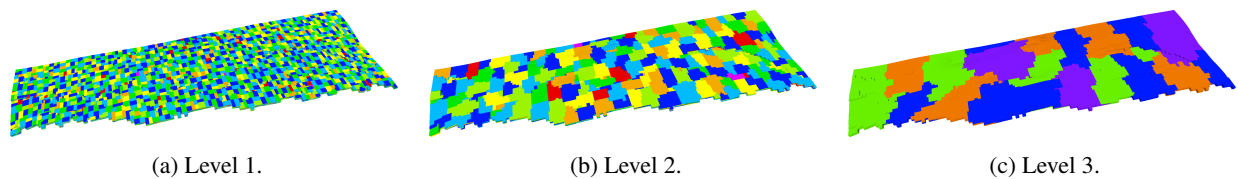

 \small
 \centering
 \begin{subfigure}[b]{0.33\textwidth}
   \centering    
   \if \generateTikzFigures 1
     \include{./pics/paper_pics/SAIGUP_lvl1}   
   \else
      \includegraphics[scale=1]{main-figure\theFigureCounter.pdf}
      \stepcounter{FigureCounter} 
   \fi
   \caption{Level 1.}
 \end{subfigure}
 \hfill
 \begin{subfigure}[b]{0.33\textwidth}
   \centering    
   \if \generateTikzFigures 1
     \include{./pics/paper_pics/SAIGUP_lvl2}   
   \else
      \includegraphics[scale=1]{main-figure\theFigureCounter.pdf}
      \stepcounter{FigureCounter}
   \fi
   \caption{Level 2.}
 \end{subfigure}
 \hfill
 \begin{subfigure}[b]{0.33\textwidth}
   \centering    
   \if \generateTikzFigures 1
     \include{./pics/paper_pics/SAIGUP_lvl3}   
   \else
      \includegraphics[scale=1]{main-figure\theFigureCounter.pdf}
      \stepcounter{FigureCounter} 
   \fi
   \caption{Level 3.}
 \end{subfigure}
\caption{Hierarchical aggregation of the SAIGUP model by METIS.}
\label{fig:saigup-partition}
\end{figure}
\begin{figure}[h!]
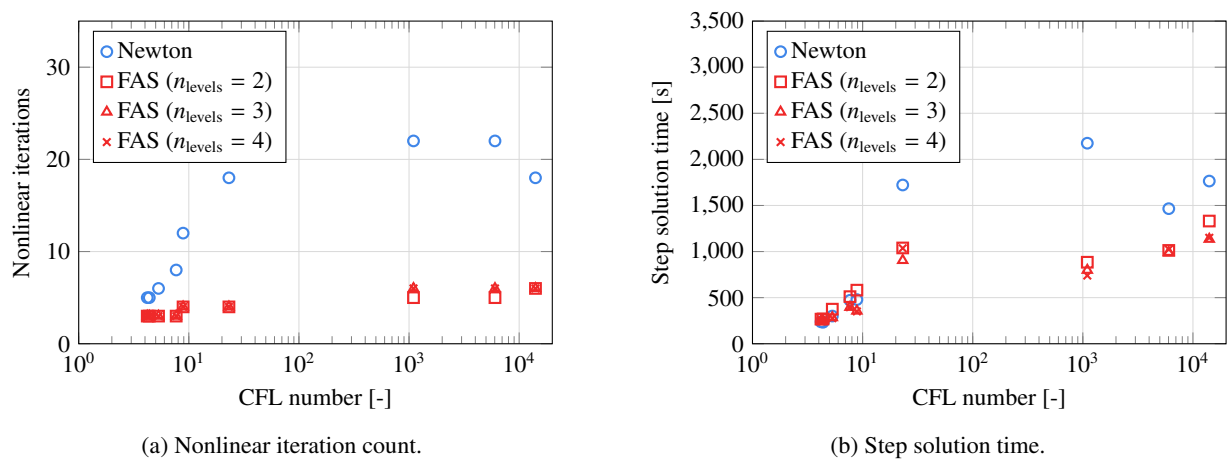

  \small
  \centering
  \begin{subfigure}[b]{0.475\textwidth}
    \centering
    \if \generateTikzFigures 1
      \include{./pics/paper_pics/SAIGUP_NLIter}    
    \else
      \includegraphics[scale=1]{main-figure\theFigureCounter.pdf}
      \stepcounter{FigureCounter}
    \fi
    \caption{Nonlinear iteration count.}
  \end{subfigure}
  \hfill
  \begin{subfigure}[b]{0.475\textwidth}
    \centering
    \if \generateTikzFigures 1
      \include{./pics/paper_pics/SAIGUP_tSolv}  
    \else
      \includegraphics[scale=1]{main-figure\theFigureCounter.pdf}
      \stepcounter{FigureCounter}
    \fi
    \caption{Step solution time.}
  \end{subfigure}
  \caption{Number of nonlinear iterations per time step and step solution time [s] as a function of CFL number for the refined SAIGUP model. FAS relies on a coarsening factor $\beta = 32$.}
  \label{fig:saigup-ml}
  
\end{figure}

We use this large test case to study the impact of the number of levels on the nonlinear behavior and total solution of FAS.
We consider up to three coarse levels, illustrated in Fig.~\ref{fig:saigup-partition}, with a coarsening factor $\beta = 32$.
The results are compared with the performance of single-level Newton for each time step in Fig.~\ref{fig:saigup-ml} and for the full simulation in Table~\ref{tab:saigup-ml-summary}.
We observe that for the three configurations considered in the table, FAS exhibits a smaller solution time than single-level Newton.
%
When we increase the number of levels from two to three, we note that the FAS nonlinear behavior slightly deteriorates, but that the FAS solution time decreases as work is shifted to less expensive computations on coarser levels.
In the most efficient multilevel configuration, the three-level FAS achieves a reduction in solution time by 39\% compared to single-level Newton.

\begin{table}[htbp]
	\caption{Solution time, $T_{\text{sol}}$ [s], and average number of nonlinear iterations per time step ($n_{\text{it}}$) for the refined SAIGUP model. FAS relies on a coarsening factor $\beta = 32$.}
	\label{tab:saigup-ml-summary}
        \small
        \centering
        \begin{tabular}{lcccc}
          \toprule
          &
          & \multicolumn{3}{c}{FAS}\\\cmidrule{3-5}
          & Newton
          & $n_{\text{levels}}=2$
          & $n_{\text{levels}}=3$
          & $n_{\text{levels}}=4$ \\\midrule
          $n_{\text{it}}$
          & 12.89 & 4.00 & 4.22 & 4.22  \\
          $T_{\text{sol}}$
          & 8,849 & 6,270 & 5,370 & 5,478 \\
          \bottomrule          
        \end{tabular}
\end{table}

\section{Concluding remarks} \label{sec:conclusion}

A nonlinear multigrid solver for two-phase flow and transport problem in a mixed fractional-flow formulation is developed.
In this formulation, the primary unknowns are the total flux, pressure, and wetting-phase saturation.
The coarse space for flux is the lowest order coarse space used in \cite{fas-spectral-diffusion}, while the coarse spaces for pressure and saturation are piecewise constant functions.
With this choice of coarse spaces, the coarse problems can be assembled with a complexity proportional to the number of cells and faces on the coarse levels, which is crucial to  arithmetic scalability.
Our numerical results show that the proposed multigrid solver exhibits a more robust nonlinear behavior than the standard single-level Newton and reduces the step solution time, especially for large CFL numbers.
This is an encouraging step to reduce the computational cost of practical large-scale reservoir simulation studies, and ensure that the time step size can be chosen based on accuracy considerations only.

Although the discussion in the current paper is based on TPFA, we remark that the proposed multigrid solver is also applicable if mimetic finite difference or mixed finite element methods are used in the discretization of \eqref{eq:model}. The resulting discrete problems will have a similar structure except that $\Mat{M}(\Vec{s})$ is no longer a diagonal matrix, cf. \cite[Chapter~6]{lie19}.
Moreover, a natural extension of the proposed solver is to use higher-order coarse spaces from \cite{fas-spectral-diffusion}, which will be explored in future work.

\section*{Acknowledgements}
\label{sec::acknow}
Funding was provided by TotalEnergies through the FC-MAELSTROM project.
Portions of this work were performed under the auspices of the U.S. Department of
Energy by Lawrence Livermore National Laboratory under Contract DE-AC52-07-NA27344 (LLNL-JRNL-826461).

\appendix

\section{The Jacobian matrix in Newton iterations}\label{app:jacobian}

Here, we provide the details of the formation of the Jacobian matrix in \eqref{eq:jacobian}.
For $\Vec{r}_\sigma^{m, \ell}$, since the nonlinear component $\Mat{M}(\Vec{s})$ has the same structure as in \cite{fas-spectral-diffusion},
$\frac{\partial \Vec{r}_\sigma^{m, \ell}}{\partial \Vec{\sigma}}$ and $\frac{\partial \Vec{r}_\sigma^{m, \ell}}{\partial \Vec{s}}$ can be formed using the approach given in Section~3.3.2 of \cite{fas-spectral-diffusion}.
Next, to find out the partial derivatives of $\Vec{r}_\sigma^{m, \ell}(\Vec{\sigma}, \Vec{s}) = T^{m, \ell}(\Vec{\sigma}, \Vec{s}) - \Vec{h}^{m, \ell}$, we can exploit the structure of $T^{m, \ell}(\Vec{\sigma}, \Vec{s})$ given in \eqref{eq:T_structure}.
Note that $T^{m, \ell}(\Vec{\sigma}, \Vec{s})$ is not differentiable due to the upwind flux.
Therefore, we consider a slightly different operator $T^{m, \ell, k}$ where the upwind direction is determined by the flux solution $\Vec{\sigma}^{\ell, k-1}$ from the previous Newton iteration:
\begin{equation}
T^{m, \ell, k}(\Vec{\sigma}, \Vec{s}) := (\Delta t_m)^{-1} W^\ell \Vec{s} + D^\ell \diag{\Vec{\sigma}}U^\ell(\Vec{\sigma}^{\ell, k-1}) f_w(\Vec{s}).
\label{eq:T_approximate}
\end{equation}
It is easy to see that
\begin{equation}
\left. \frac{\partial \Vec{r}_s^{m, \ell}}{\partial \Vec{s}} \right|_{(\Vec{\sigma}, \Vec{s}) = (\Vec{\sigma}^{\ell, k-1}, \Vec{s}^{\ell, k-1})}
\approx \left. \frac{\partial T^{m, \ell, k}}{\partial \Vec{s}} \right|_{(\Vec{\sigma}, \Vec{s}) = (\Vec{\sigma}^{\ell, k-1}, \Vec{s}^{\ell, k-1})}
= (\Delta t_m)^{-1} W^\ell + D^\ell \diag{\Vec{\sigma}^{\ell, k-1}}U^\ell(\Vec{\sigma}^{\ell, k-1}) \left. \frac{d \left( f_w(\Vec{s} ) \right)}{d\Vec{s}} \right|_{\Vec{s} = \Vec{s}^{\ell, k-1}} .
\end{equation}
To obtain $\frac{\partial \Vec{r}_s^{m, \ell}}{\partial \Vec{\sigma}} $, we first note that $T^{m, \ell, k}$ in \eqref{eq:T_approximate} can be rearranged to be
\begin{equation}
T^{m, \ell, k}(\Vec{\sigma}, \Vec{s}) = (\Delta t_m)^{-1} W^\ell \Vec{s} + D^\ell \diag{U^\ell(\Vec{\sigma}^{\ell, k-1}) f_w(\Vec{s})}\Vec{\sigma}.
\label{eq:T_approximate_rearranged}
\end{equation}
Therefore,
\begin{equation}
\left. \frac{\partial \Vec{r}_s^{m, \ell}}{\partial \Vec{\sigma}} \right|_{(\Vec{\sigma}, \Vec{s}) = (\Vec{\sigma}^{\ell, k-1}, \Vec{s}^{\ell, k-1})}
\approx \left. \frac{\partial T^{m, \ell, k}}{\partial \Vec{\sigma}} \right|_{(\Vec{\sigma}, \Vec{s}) = (\Vec{\sigma}^{\ell, k-1}, \Vec{s}^{\ell, k-1})}
= D^\ell \diag{U^\ell(\Vec{\sigma}^{\ell, k-1}) f_w(\Vec{s}^{\ell, k-1})}.
\end{equation}



\end{document}